\documentclass[11pt]{article}

\usepackage[margin=.7in]{geometry}
\usepackage{amsfonts}
\usepackage{amsmath}
\usepackage{amsthm}

\usepackage{hyperref}

\usepackage{xcolor}

\newcommand{\mi}{\textrm i}

\newtheorem{thm}{Theorem}

\newtheorem{prop}{Proposition}
\newtheorem{dfn}{Definition}
\newtheorem{rmk}{Remark}
\newtheorem{lmm}{Lemma}

\usepackage{./algorithmicx}
\usepackage[linesnumbered,ruled,vlined]{./algorithm2e}
\usepackage[noend]{./algpseudocode}

\usepackage{placeins}

\title{Local Taylor-based polynomial quasi-Trefftz spaces for scalar linear equations
}

\author{Lise-Marie Imbert-G\'erard}

\begin{document}
\maketitle
\tableofcontents
\newpage



{\bf Abstract}\\
Trefftz-type of Galerkin methods for numerical PDEs use discrete spaces of problem-dependent functions.
While Trefftz methods leverage discrete spaces of local exact solutions to the governing PDE, Taylor-based quasi-Trefftz methods leverage discrete spaces of local approximate solutions to the governing PDE. This notion of approximate solution, understood in the sense of a small Taylor remainder, is defined for differential operators with smooth variable coefficients.
In both cases, it is possible to use discrete spaces much smaller than standard polynomial space to achieve the same orders of approximation properties.

The present work is the first systematic study of local Taylor-based polynomial quasi-Trefftz spaces 
characterized as the kernel of the quasi-Trefftz operator, defined as the composition of Taylor truncation with the differential operator.
The proposed linear algebra framework reveals the general structure of this linear operator and applies to any non-trivial linear scalar differential operator with smooth coefficients. It results in a fully explicit procedure to construct a local quasi-Trefftz basis valid in all dimension and for operators of any order, guaranteeing a minimal computational cost for the construction of these equation-dependent bases.


The local quasi-Trefftz space is formally defined as the kernel of a linear operator between spaces of polynomials. The systematic approach relies on a detailed study of the structure of this operator, strongly leveraging the graded structure of polynomial spaces.

\section{Introduction}
Trefftz Discontinuous Galerkin (TDG) methods are a class of high-order numerical methods for Partial Differential Equation (PDE) problems, including boundary value problems and initial boundary value problems.  They are named after the original work of E. Trefftz \cite{Trefftz}.
As Discontinuous Galerkin (DG) methods, they rely on a weak formulation of the problem. 
As Trefftz methods they rely on function spaces of exact local solutions to the governing PDE.
Thanks to adequate discrete spaces of exact local solutions to the governing PDE, TDG methods -- compared to DG methods with standard polynomial discrete spaces -- can achieve the same order of convergence with much fewer degrees of freedom, see \cite{survey}. 

In the field of numerical PDEs, several methods leverage problem-dependent function spaces. This is for instance the case of the Discontinuous-Petrov-Galerkin method \cite{DPG-0,DPG-Dem} and of multiscale methods \cite{Peters-1,Peters-2}. Instead of computed problem-dependent basis functions, TDG methods rely on bases of exact solutions to the governing equation \cite{Git,survey}. While plane wave functions are convenient exact solutions to various differential equations modeling wave propagation, for many governing equations bases of exact solutions might not be available.
In this case, various flavors of Trefftz-like methods relying on discrete spaces of such approximate rather than exact solutions to the governing PDE have been introduced. The main idea is to relax the so-called Trefftz property, considering instead a local approximate version:
$$
\mathcal L \phi \approx 0
\text{ as opposed to } 
\mathcal L \phi = 0
,$$
where $\mathcal L$ is the partial differential operator of the governing PDE of interest.
For instance, the embedded Trefftz method \cite{emb} relies on a Galerkin projection onto a subspace of Trefftz type, while a quasi-Trefftz method \cite{Seb} relies on an impedance trace approximation, both leveraging local solvers for the computation of basis functions.

Another type of quasi-Trefftz methods has been introduced by the author and collaborators, interpreting the quasi-Trefftz property in the sense of a truncated Taylor expansion: the image of a function by the partial differential operator is not equal to zero, but its Taylor expansion up to a certain order is equal to zero. 
Such methods have been developed for the Helmholtz equation \cite{IGD-GPW,IGM-GPW}, for the time-dependent wave equation \cite{IGMS-time}, for a family of elliptic equations \cite{IGMPS-dar}, for the Schr\"odinger equation \cite{schrod}. 

Both for Trefftz and for quasi-Trefftz methods, the definition, construction and approximation properties of discrete spaces have also been a focus of interest, 
for Schr\"odinger equation \cite{SchrodSpace},
for Helmholtz equation \cite{evan,MYinterp},
for the convected Helmholtz equation \cite{convHelm},
and more general homogeneous scalar equations \cite{roadmap,ampbased},
as well as in a vector-valued case \cite{vect}.
These references study either polynomial or wave-like ansatz for the basis functions.

Taylor-based polynomial quasi-Trefftz functions were first introduced in \cite{IGMS-time} for the wave equation with variable sound speed;  they were studied for 3D second-order problems in \cite{convHelm}; they were further studied in \cite{IGMPS-dar}. 
Each of these works relies on an assumption on the highest order terms in the differential operator.\footnote{Either $c_{\gamma\mathbf e_1}\left( {\mathbf x}  \right) \neq 0$ for the general order $\gamma$ operator $\displaystyle \mathcal L := \sum_{m=0}^\gamma \sum_{ |\mathbf j|=m} c_{\mathbf j}\left( {\mathbf x}  \right) \partial_{{\mathbf x}}^{\mathbf j}$, or $\begin{bmatrix}  c_{2e_1} & \frac12 c_{e_1+e_2}  & \frac12 c_{e_1+e_3}  \\ \frac12 c_{e_1+e_2}  &  c_{2e_2}  & \frac12  c_{e_2+e_3}  \\
  \frac12 c_{e_1+e_3}  &\frac12  c_{e_2+e_3}   &  c_{2e_3} 
\end{bmatrix}$ non-singular for second order operators $\mathcal L=\sum_{i\in\mathbb N_0^3; |i|\leq 2} c_{\mathbf j}\left( {\mathbf x}  \right) \partial_{{\mathbf x}}^i$. }
{\bf The present work is the first systematic study of local Taylor-based polynomial quasi-Trefftz spaces 
characterized as the kernel of the quasi-Trefftz operator, defined as the composition of Taylor truncation with the differential operator.
The proposed linear algebra framework reveals the general structure of this linear operator and applies to any non-trivial linear scalar differential operator with smooth coefficients. It results in a fully explicit procedure to construct a local quasi-Trefftz basis valid in all dimension and for operators of any order.}
This explicit procedure guarantees a minimal computational cost for the construction of a global quasi-Trefftz space to implement a quasi-Trefftz method.
Moreover, even though the general tools developed here do not directly apply to quasi-Trefftz spaces for different ansatz, they still provide fundamental ideas that can be adapted beyond the polynomial setting, see \cite{GPWgen}.

In the field of numerical integral equations, polynomial solutions of linear constant-coefficient PDEs with a polynomial right-hand side can be leveraged to evaluate volume integral operators \cite{Bonnet++}.
The construction of such solutions is tackled in \cite{OtherPols,BonnetPols}. The former focuses on second order equation and proposes a procedure leveraging the constant term of the differential operator. The latter tackles operators of arbitrary orders and proposes a procedure hinging on the elementary part of the differential operator, hence generalizing the previous approach. 
There the general methodology crucially relies on graded structure of polynomial spaces.
The general framework presented here in Sections \ref{ssec:GFsurj}, \ref{ssec:GFker} and \ref{ssec:GFbas} also relies on a key additive decomposition of the linear operator into a block-diagonal part and the corresponding remainder, where the block structure follows the graded structure of polynomial spaces.
This framework is applied to quasi-Trefftz functions and variable-coefficients PDEs in Sections \ref{ssec:surj}, \ref{ssec:ker} and \ref{ssec:bas}.
 While in this case the diagonal operator is the principal part of the differential operator, the framework can also be applied with the diagonal operator being the elementary part of the differential operator to cover the case of exact polynomial solutions (for constant coefficient operators and polynomial right-hand sides).

\subsection{Notation}
Throughout the article, 
while $\mathbb N$ refers to the set of positive integers, $\mathbb N_0$ refers to the set of non-negative integers;
for any $(x,y)\in\mathbb R$ with $x\leq y$ then $[\![x,y]\!] \subset \mathbb Z$ is the set of integers between $x$ and $y$;
the number of elements in a set $E$ is denoted $\dim  E$;
 $\mathcal K$ denotes the kernel and $\mathcal R$ range.
Moreover, 
$g$ refers to an integer,  
$d$ refers to the dimension of the ambient space, with $\mathbf x=(x_g,g\text{ from } 1 \text{ to } d)\in\mathbb R^d$, 
$\mathbf i=(i_g,g\text{ from } 1 \text{ to } d)\in(\mathbb N_0)^d$ refers to a multi-index;
partial derivatives are denoted $\partial^{\mathbf i}=\partial_{x_1}^{i_1}\cdots\partial_{x_d}^{i_d} $, 
and finally the canonical basis of $\mathbb R^d$ is denoted $\{\mathbf e_g, 1\leq g\leq d\}$. The relations $\leq$ and $<$ between multi-indices is interpreted to hold component-wise.

The space of polynomials of degree at most equal to $p$ in $d$ variables on the field $\mathbb K=\mathbb R$ or $\mathbb C$ is denoted $\mathbb P_p$.
For any $l\in\mathbb N_0$, the space of homogeneous polynomials of degree $l$ is denoted $\widetilde{\mathbb P}_l$. They are related by $\displaystyle \mathbb P_p = \bigoplus_{l=0}^p \widetilde{\mathbb P}_l$.  
The canonical monomial is denoted $\displaystyle \mathbf X^{\mathbf i}=\prod_{g=1}^d (X_{g})^{i_g}$ for any $\mathbf i\in(\mathbb N_0)^d$.

The work presented here is limited to {\bf local} properties. 
More precisely, we only consider spaces defined in the neighborhood of a point $\mathbf x_0\in\mathbb R^d$, so $\mathbf x_0$ will be omitted in the following notation. 
For $n\in\mathbb N_0$, the space of functions with continuous derivatives up to order $n$ in a neighborhood of $\mathbf x_0$ is denoted $\mathcal C^n$.
For $k\in\mathbb N_0$, for a $\varphi\in\mathcal C^k$ function $\varphi$, its $k$th Taylor polynomial at $\mathbf x_0$, denoted $T_k[\varphi] \in\mathbb P_k$, is given by
$$
T_k[\varphi] =\sum_{l=0}^k \sum_{|\mathbf i| = l } \frac 1{\mathbf i!} \partial^{\mathbf i} \varphi(\mathbf x_0) (\mathbf X-\mathbf x_0)^{\mathbf i}.
$$

\subsection{
Quasi-Trefftz spaces}
First, the precise notion of Taylor-based quasi-Trefftz property of a function, as first introduced in \cite{IGD-GPW}, is reminded here. This is of course a local property as it is expressed in terms of a Taylor expansion.
\begin{dfn}
Assume $q\in\mathbb N_0$.
Given a linear partial differential operator $\mathcal L$ with smooth variable coefficients in $\mathcal C^{q}$, the Taylor-based quasi-Trefftz property of order $q$ for a function $\varphi$ is
$$
T_q[\mathcal L  \varphi] = 0_{\mathbb P_q}.
$$
\end{dfn}
A polynomial quasi-Trefftz space can then be defined as a space of polynomials whose image by the differential operator satisfies an adequate quasi-Trefftz property. 
\begin{dfn}
Assume $\gamma\in\mathbb N$ and $p\in\mathbb N$ with $p\geq \gamma$.
For a linear partial differential operator $\mathcal L$ of order $\gamma$ with smooth variable coefficients in $\mathcal C^{p-\gamma}$, the polynomial quasi-Trefftz space of order $p$ is defined as the space of polynomials of degree at most equal to $p$ whose image by the differential operator $\mathcal L$ satisfies the quasi-Trefftz property for $q=p-\gamma$. This quasi-Trefftz space can be written as
$$
\mathbb Q\mathbb T_p = \{\Pi\in \mathbb P_p | T_{p-\gamma}[\mathcal L  \Pi] = 0_{\mathbb P_{p-\gamma}} \},
$$
 it is the kernel of the so-called quasi-Trefftz operator
\begin{equation}\label{eq:op}
\begin{array}{rll}
\mathcal D_p:& \mathbb P_p& \to  \mathbb P_{p-\gamma}\\
& \Pi&\mapsto T_{p-\gamma}[\mathcal L  \Pi] 
\end{array}
\end{equation}
\end{dfn}
\begin{rmk}\label{rmk:notDO}
While $\mathcal L$ is a differential operator, in general it is not the case of the quasi-Trefftz operator $\mathcal D_p$. 
\end{rmk}

\begin{rmk}
Unlike standard polynomial spaces, the quasi-Trefftz spaces are not embedded,
$$
\forall p>q,\mathbb{QT}_q \not\subset \mathbb{QT}_p,
$$
since for $\Pi\in\mathbb{QT}_q$ there is no condition imposed on $T_p \Pi-T_q \Pi$,
 however
$$
\forall p>q,\forall \Pi\in\mathbb{QT}_p,  T_q\Pi\in\mathbb{QT}_q.
$$
\end{rmk}
In the previous definition, the order of the quasi-Trefftz property in terms of the polynomial degree $p$ and the differential operator's order $\gamma$ is chosen to be $p-\gamma$. This can be explained as follows.
Given a variable-coefficient partial differential operator defined as
$\displaystyle \mathcal L_c := \sum_{m=0}^\gamma \sum_{ |\mathbf j|=m} c_{\mathbf j}\left( {\mathbf x}  \right) \partial_{{\mathbf x}}^{\mathbf j}$ 
with variable coefficients
$\displaystyle c_{\mathbf j} = \sum_{n=0}^{p-\gamma} \sum_{|\mathbf k| = n } c_{\mathbf j,\mathbf k} (\mathbf X-\mathbf x_0)^{\mathbf k} + O(\| \mathbf X-\mathbf x_0 \|^{p-\gamma+1})$. 
For all polynomial $\Lambda\in\mathbb P_p$ 
with 
$\displaystyle \Lambda =  \sum_{l=0}^{p} \sum_{|\mathbf i| = l } \lambda_{\mathbf i} (\mathbf X-\mathbf x_0)^{\mathbf i}
$,
for $0\leq m \leq \gamma, |\mathbf j|=m$, then
$$
 c_{\mathbf j}\left( {\mathbf x}  \right) \partial_{{\mathbf x}}^{\mathbf j}\sum_{l=0}^{p} \sum_{|\mathbf i| = l } \lambda_{\mathbf i} (\mathbf X-\mathbf x_0)^{\mathbf i}
=
c_{\mathbf j}\left( {\mathbf x}  \right)
\left(  
\sum_{l=0}^{p-m}
 \sum_{|\mathbf i| = l } \frac{(\mathbf i+\mathbf j)!}{\mathbf i!} \lambda_{\mathbf i+\mathbf j} (\mathbf X-\mathbf x_0)^{\mathbf i}  \right)
$$
\begin{align*}
c_{\mathbf j}\left( {\mathbf x}  \right) \partial_{{\mathbf x}}^{\mathbf j}\sum_{l=0}^{p} \sum_{|\mathbf i| = l } \lambda_{\mathbf i} (\mathbf X-\mathbf x_0)^{\mathbf i}
=
\sum_{l=0}^{p-m}  \sum_{|\mathbf i|=l}
 \left( \sum_{n= 0}^{l} \sum_{|\mathbf k|=n}    c_{\mathbf j,\mathbf k} \frac{(\mathbf i-\mathbf k+\mathbf j)!}{(\mathbf i-\mathbf k)!} \lambda_{\mathbf i-\mathbf k+\mathbf j} \right)
(\mathbf X-\mathbf x_0)^{\mathbf{i}} \\
 +
 \sum_{l=p-m+1}^\infty  \sum_{|\mathbf i|=l}
 \left( \sum_{n=l-(p-m) }^{l} \sum_{|\mathbf k|=n}    c_{\mathbf j,\mathbf k}  \frac{(\mathbf i-\mathbf k+\mathbf j)!}{(\mathbf i-\mathbf k)!} \lambda_{\mathbf i-\mathbf k+\mathbf j} \right)
(\mathbf X-\mathbf x_0)^{\mathbf{i}} 
\end{align*}
In particular, in the Taylor expansion of $\mathcal L_c \Lambda$, each of the terms of order $m$ from 0 to $p-\gamma$ includes contributions to their coefficient of the form
$$
 \left( \sum_{n= 0}^{l} \sum_{|\mathbf k|=n}    c_{\mathbf j,\mathbf k} \frac{(\mathbf i-\mathbf k+\mathbf j)!}{(\mathbf i-\mathbf k)!} \lambda_{\mathbf i-\mathbf k+\mathbf j} \right),
 \quad \forall \mathbf j\in\mathbb N_0^d \text{ with } |\mathbf j|=\gamma,
$$
while the terms of order $m$ at least equal to $p-\gamma+1$ do not. 
Similarly,  in the Taylor expansion of $\mathcal L_c \Lambda$, the terms of order $m$ from 0 to $p-m$ include contributions to their coefficient of the form
$$
 \left( \sum_{n= 0}^{l} \sum_{|\mathbf k|=n}    c_{\mathbf j,\mathbf k} \frac{(\mathbf i-\mathbf k+\mathbf j)!}{(\mathbf i-\mathbf k)!} \lambda_{\mathbf i-\mathbf k+\mathbf j} \right),
 \quad \forall  \mathbf j\in\mathbb N_0^d \text{ with } |\mathbf j|<\gamma,
$$
Hence, in particular, they contribute to the terms of order $m$ from 0 to $p-\gamma$.
Therefore,  in the Taylor expansion of $\mathcal L_c \Lambda$, only the terms of order $m$ from 0 to $p-\gamma$ contain the same contributions, while terms of order $m$ larger than $p-\gamma$ do not.
This motivates the choice to define the quasi-Trefftz space as the kernel of $T_{p-\gamma}\mathcal L$, in order for all the equations defining the quasi-Trefftz space to share the same structure.


\subsection{Best approximation properties}
Best approximation properties of discrete spaces are key in the proof of convergence of Galerkin-type methods. Here, following the scope of this article, we limit our attention to local properties, that would then be applied element per element in the Galerkin setting. 

In the case of polynomial quasi-Trefftz spaces, the best approximation results entirely rely on the following property of exact solutions to the homogeneous PDE, namely that their Taylor polynomial belongs to the quasi-Trefftz space, which does not require specific hypothesis on the partial differential operator. 
This kind of approximation results for quasi-Trefftz spaces were first proved for Generalized Plane Wave spaces \cite{MYinterp}. 
\begin{prop}\label{prop:TuinQT}
Assume $\gamma\in\mathbb N$ and $p\in\mathbb N$ with $p\geq \gamma$.
Consider a linear partial differential operator $\mathcal L$ of order $\gamma$ with smooth variable coefficients in $\mathcal C^{p-\gamma}$.
For any function $u\in\mathcal C^{p+\gamma}$ such that $\mathcal L u=0$ in a neighborhood of $\mathbf x_0$, the Taylor polynomial of order $p$ of $u$ belongs to the quasi-Trefftz space of order $p$: $T_{p}u\in\mathbb{QT}_p$. 
\end{prop}
\begin{rmk}
Of course, in general, $T_p \partial^{\mathbf i} \Pi\neq  \partial^{\mathbf i} T_p\Pi$ for polynomials $\Pi$. 
This can be illustrated by the fact that $ \partial^{\mathbf i} T_p\Pi \in\mathbb P_{p-|\mathbf i|}$ while $T_p \partial^{\mathbf i} \Pi\in\mathbb P_p$, hence as long as $\deg (T_p \partial^{\mathbf i} \Pi) >p-|\mathbf i|$ then clearly $T_p \partial^{\mathbf i} \Pi\neq  \partial^{\mathbf i} T_p\Pi$.
For instance, for $p>0$ consider $\Pi=\mathbf X^{(p+1)\mathbf e_1}$. In this case $T_p\Pi=0$ so $ \partial^{\mathbf i} T_p\Pi=0$ for any $\mathbf i$, while for $\mathbf i=2\mathbf e_1$ in this case $\partial^{\mathbf i} \Pi\in\mathbb P_{p-2}$ so $T_p \partial^{\mathbf i} \Pi= \partial^{\mathbf i} \Pi$, and $\partial^{\mathbf i} \Pi=(p+1)p\mathbf X^{(p-1)\mathbf e_1}$, which is not equal to zero.
%
\end{rmk}
\begin{proof}
Since $u\in\mathcal C^{p-\gamma}$, then $\mathcal L u=\mathcal LT_{p}u+\mathcal L(u-T_p u)$, but $\mathcal LT_{p}u\in\mathbb P_{p-2}$ while $\mathcal L(u-T_p u)\in\mathcal K (T_{p-2})$. 
As a consequence, $T_{p-2}\mathcal L u=T_{p-2}\mathcal LT_{p}u$. 
Since moreover $\mathcal Lu=0$,  then $T_{p-2}\mathcal L(T_{p} u)=0$. 
This concludes the proof as it shows that $T_{p} u\in\mathbb{QT}_p$.
\end{proof}
By contrast, the choice of the Generalized Plane Wave (GPW) ansatz together with the quasi-Trefftz property does not define a finite dimensional space. Instead of Proposition \ref{prop:TuinQT}, the adequate Taylor polynomial of exact solutions to the homogeneous PDE can be shown to belong to a space of Taylor polynomials of functions satisfying the quasi-Trefftz property, independently of the ansatz. 
Both in the polynomial and the GPW cases, this property of a Taylor polynomial of exact solutions to the homogeneous PDE is essential to best approximation properties.

The fundamental best approximation property of polynomial quasi-Trefftz spaces follows. 
\begin{prop}
Assume $d\in\mathbb N$, $\mathbf x_0\in\mathbb R^d$, $\gamma\in\mathbb N$ and $p\in\mathbb N$ with $p\geq \gamma$.
Consider a linear partial differential operator $\mathcal L$ of order $\gamma$ with smooth variable coefficients in $\mathcal C^{p-\gamma}$.
For any function $u\in\mathcal C^{p+\gamma}$ such that $\mathcal L u=0$ in a neighborhood of $\mathbf x_0$, there exists a quasi-Trefftz function $u_a\in\mathbb{QT}_p$ such that, in a neighborhood of $\mathbf x_0$, 
$$
\left\{\begin{array}{l}
|u(\mathbf x)- u_a(\mathbf x)|\leq C |\mathbf x-\mathbf x_0|^{p+1}\\
|\nabla u(\mathbf x)- \nabla u_a(\mathbf x)|\leq C |\mathbf x-\mathbf x_0|^{p}
\end{array}\right.
$$
for some constant $C$ independent of $\mathbf x$, but depending on $u$, $\mathbf x_0$ and $p$.
\end{prop}
\begin{proof}
According to Proposition \ref{prop:TuinQT}, $T_p u\in\mathbb{QT}_p$.
The result follows with $u_a=T_p u$.
\end{proof}

In \cite{convHelm}, a single proof of these approximation properties in dimension $d=3$ was presented for polynomial quasi-Trefftz spaces and for two types of GPW spaces. The backbone of that proof is a space of Taylor polynomials of functions satisfying the quasi-Trefftz property, defined independently of the ansatz. Unlike here, that proof relies on hypotheses on the differential operator $\mathcal L$; one hypothesis is a particular case of the central hypothesis on which the rest of this work relies, while the other one is required only to handle the GPW case.

As a final remark on quasi-Trefftz approximation properties, it is important to highlight how they compare to approximation properties of standard polynomial spaces. 
On the one hand, the space of functions that can be approximated are different: quasi-Trefftz spaces $ \mathbb Q\mathbb T_n$ can only approximate exact solutions to the governing equation while polynomial spaces $ \mathbb P_n$ can approximate functions based on their smoothness. This is only natural, since quasi-Trefftz spaces are tailored to a single equation, but it is independent of the desired order of approximation.
On the other hand, $ \mathbb Q\mathbb T_n\subset\mathbb P_n$ and the difference between their dimensions depends on the order $\gamma$ of the differential operator: 
$$
\dim  \mathbb Q\mathbb T_n =\dim  \mathbb P_n -\dim  \mathbb P_{n-\gamma} 
$$
so $\dim  \mathbb Q\mathbb T_n =O(n^{d-1})$ while $\dim  \mathbb P_n =O(n^{d})$, see examples in Table \ref{tab:dims}. 

\begin{table}[h!]
  \begin{center}
    \caption{Comparing dimensions of discrete spaces with the same order $n$ of approximation properties in dimension $d$, between the polynomial space $\mathbb P_n$ and some quasi-Trefftz space $ \mathbb Q\mathbb T_n$.}\label{tab:dims}
    \label{tab:table1}
{\renewcommand{\arraystretch}{1.75}%
    \begin{tabular}{c|ccc} 
   &   \textbf{$d=2$} & \textbf{$d=3$} & \textbf{$d$} \\
      \hline
   $ \mathbb P_n$&  $\displaystyle\frac{(n+1)(n+2)}{2}=O\left(n^2\right)$ & $\displaystyle\frac{(n+1)(n+2)(n+3)}{6}=O\left(n^3\right)$ & $O(n^{d})$ \\
   $ \mathbb Q\mathbb T_n$ with $\gamma =2$&  $2n+1=O\left(n\right)$ & $(n+1)^2=O\left(n^2\right)$& $O(n^{d-1})$  \\
   $ \mathbb Q\mathbb T_n$ with $\gamma =3$&  $3n=O\left(n\right)$ & $\frac{3n^2+3n+2}2=O\left(n^2\right)$& $O(n^{d-1})$ 
    \end{tabular}}
  \end{center}
\end{table}

\subsection{Highlights}

The PDEs of interest are scalar and homogeneous.
The proposed notion of quasi-Trefftz function relies on a Taylor expansion: the image of a quasi-Trefftz function by the differential operator is not zero, but its Taylor expansion is zero up to a certain order.
Hence all the results in this work are local, around a point where the PDE coefficients are sufficiently smooth. 
The PDEs of interest are also linear, leading to a linear quasi-Trefftz operator at the center of this study.

The final goal of this article is to propose a generic algorithm to construct polynomial quasi-Trefftz bases  independently of the dimension $d$, see Section \ref{sec:bas}. This is possible thanks to a comprehensive study of the fundamental structure of polynomial quasi-Trefftz spaces, including two main aspects:
\begin{itemize}
\item the surjectivity of quasi-Trefftz operator, see Section \ref{sec:surj};
\item the characterization of quasi-Trefftz space as kernel of the quasi-Trefftz operator, see Section \ref{sec:ker}.
\end{itemize}
Relying on the graded structure of polynomial spaces, the theoretical foundation for all results is presented  in the framework of block linear algebra, and more particularly to the underlying block-triangular structure of relevant operators. 

These three sections rely on a central hypothesis for the differential operator. 
Given a differential operator $\displaystyle \mathcal L := \sum_{m=0}^\gamma \sum_{ |\mathbf j|=m} c_{\mathbf j}\left( {\mathbf x}  \right) \partial_{{\mathbf x}}^{\mathbf j}$, this central hypothesis states that the operator $\displaystyle \mathcal L_*:=  \sum_{ |\mathbf j|=\gamma} c_{\mathbf j}\left( {\mathbf x_0}  \right) \partial_{{\mathbf x}}^{\mathbf j}$  
defines a surjective operator on polynomials spaces. 
 Section \ref{sec:apps} focuses on this question. It shows that all non-trivial constant-coefficient homogeneous differential operators satisfy this property, and then concludes with Theorem \ref{thm:qTccl} summarizing how to construct a quasi-Trefftz basis for any scalar homogeneous linear smooth-coefficient PDE of order exactly $\gamma\in\mathbb N$ at $\mathbf x_0$ in any dimension $d$.

\section{Quasi-Trefftz operator
}\label{sec:surj}
Since the polynomial quasi-Trefftz space is introduced as the kernel of a linear operator, let us focus first on thisquasi-Trefftz operator $\mathcal D_p$ defined in \eqref{eq:op}.
The goal of this section is to prove that $\mathcal D_p$ is surjective, and to provide a formula for a right inverse operator.

Examining the operator $\mathcal D_P$ from the point of view of the graded structure of polynomial spaces is key.
Since for any $\{\Pi_n\in \widetilde{\mathbb P}_n, n\in[\![0,p]\!]\}$ 
$$
\mathcal D_p \left(\sum_{n = 0}^p\Pi_n\right)=\sum_{n = 0}^p\mathcal D_p ( \Pi_n)
\text{ and }
 \mathcal D_p(\Pi_n) = \sum_{m=\max(0,n-(p-\gamma))}^{\min(n,\gamma)} \sum_{ |\mathbf j|=m}  \sum_{\ell=n-m}^{p-\gamma} \underbrace{c_{\mathbf j,\ell-(n-m)} \partial_{\mathbf x}^{\mathbf j}\Pi_n}_{\in\widetilde{\mathbb P}_\ell} \text{ (see Appendix \ref{app:newexpl})},
$$
 then for any $n\in[\![\gamma,p]\!]$ the range of $\mathcal D_p|_{\widetilde{\mathbb P}_{n}}$ is a subset of $\bigoplus_{k=n-\gamma}^{p-\gamma} \widetilde{\mathbb P}_k$
while the range of $\mathcal D_p$ is a subset of $\mathbb P_{p-\gamma}$. 
This in particular implies that the operator $\mathcal D_p$ restricted to $\bigoplus_{\gamma}^p\widetilde{\mathbb P}_{\ell}$ has a block-triangular structure.
As such, its surjectivity is equivalent to the surjectivity of its diagonal blocks, which precisely correspond to the operator $\mathcal L_*$, and a forward substitution procedure paves the way to a right inverse operator.

\subsection{General framework 
}\label{ssec:GFsurj}
This section introduces the general framework of block-triangular linear operators that captures the essential properties of the quasi-Trefftz operator.
\begin{thm}\label{thm:surj}
Consider 
$s\in\mathbb N$ as well as a linear 
  operator $\mathcal T:\mathbb A\to \mathbb B$ between
 two graded vector spaces 
\begin{equation}
\mathbb A = \bigoplus_{n=0}^{s} \mathbb A_n\text{ and }\mathbb B = \bigoplus_{l=0}^{s} \mathbb B_l.
\end{equation}
Assuming that the operator can be decomposed as the sum of two operators $\mathcal T_*$ and $\mathcal R$ with
\begin{enumerate}
\item
\begin{enumerate}
\item\label{rangeT*}
for all $n\in[\![0,s]\!]$, $\mathcal T_*(\mathbb A_n)\subset \mathbb B_n$;
\item $\mathcal T_*$ is surjective;
\end{enumerate}
\item 
\begin{enumerate}
\item\label{rangeR}
for all $n\in[\![0,s-1]\!]$, $\mathcal R(\mathbb A_n)\subset \bigoplus_{l=n+1}^{s} \mathbb B_l$;
\item\label{decompR}
$\mathcal R(\mathbb A_s) = \{0_{\mathbb B}\}$;
\end{enumerate}
\end{enumerate}
then $\mathcal T $
  is surjective.
  
  Moreover, if $\mathcal S_*$ denotes a right inverse of $\mathcal T_*:\mathbb A\to \mathbb B$, such that for all $n\in[\![0,s]\!]$, $\mathcal S_*(\mathbb B_n)\subset \mathbb A_n$, then a right inverse of $\mathcal T$ is the operator $\mathcal S:\mathbb A\to \mathbb B$ given by $\mathcal S = \mathcal S_* \circ\sum_{l=0}^{s}[- \mathcal R \circ\mathcal S_*]^l$. 
\end{thm}
The graded nature of the vector spaces $A$ and $B$ naturally corresponds, in appropriate bases, to a block structure of the operator $\mathcal T$'s matrix.
Hence the two hypotheses can naturally be interpreted in terms of blocks of such a matrix, with column and row blocks numbered respectively with increasing $n$ and $l$: they correspond to a block-triangular structure, with diagonal blocks defining the operator $\mathcal T_*$ and strictly lower-triangular blocks defining the operator $\mathcal R$.
The proof of this theorem will then naturally follow a forward substitution procedure.
\begin{rmk}
Any operator $\mathcal T$ satisfying hypotheses \ref{rangeT*}, \ref{rangeR} and \ref{decompR} is surjective if and only if $\mathcal T_*$ is surjective.
\end{rmk}
\begin{proof}
As a preliminary, for all $n\in[\![0,s-1]\!]$, the operator $\mathcal R|_{\mathbb A_n}:\mathbb A_n\to \bigoplus_{l=n+1}^{s} \mathbb B_l$ can be decomposed as the sum of $s-n$ operators $\mathcal R_{ln}:\mathbb A_n\to \mathbb B_l$ for all $l\in[\![n+1,s]\!]$:
$$
\forall x_n\in \mathbb A_n, \mathcal R(x_n) = \sum_{l=n+1}^s \mathcal R_{ln}(x_n).
$$
Then the operator $\mathcal T$ is surjective if and only if for any $\{y_l\in \mathbb B_l, l\in[\![0,s]\!]\}$ there exists $\{x_n\in \mathbb A_n, n\in[\![0,s]\!]\}$ such that
\begin{equation*}
 \mathcal T\left(\sum_{n=0}^s x_n\right)=\sum_{l=0}^s y_l 
\Leftrightarrow\sum_{n=0}^s \mathcal T_*\left(x_n\right)+  \sum_{n=0}^{s-1} \mathcal R\left(x_n\right) = \sum_{l=0}^s y_l .
\Leftrightarrow\sum_{n=0}^s \mathcal T_*\left(x_n\right)+  \sum_{n=0}^{s-1} \sum_{l=n+1}^s \mathcal R_{ln}(x_n) = \sum_{l=0}^s y_l
,
\end{equation*}
$$
\Leftrightarrow
\sum_{l=0}^s \underbrace{\mathcal T_*\left(x_l\right)}_{\in \mathbb B_l}+  \sum_{l=1}^s \sum_{k=0}^{l-1} \underbrace{\mathcal R_{lk}(x_k)}_{\in \mathbb B_l} = \sum_{l=0}^s \underbrace{y_l}_{\in \mathbb B_l}
\Leftrightarrow
\left\{
\begin{array}{l}
\mathcal T_*(x_0) = y_0, \\
\displaystyle 
 \mathcal T_*\left(x_l\right)= y_l-\sum_{k=0}^{l-1} \mathcal R_{lk}(x_k), \ \forall l\in[\![1,s]\!].
\end{array}
\right.
$$
This last system can be solved for values of $n$ increasing from $0$ to $s$ since $\mathcal T_*$ is surjective and the right-hand side of each equation is known from $y$ and from previous equations. 
Hence for any $\{y_l\in \mathbb B_l, l\in[\![0,s]\!]\}$ there exists indeed $\{x_n\in \mathbb A_n, n\in[\![0,s]\!]\}$ such that $ \mathcal T\left(\sum_{n=0}^s x_n\right)=\sum_{l=0}^s y_l$: $\mathcal T$ is surjective.

Given a right inverse of $\mathcal T_*$ denoted $\mathcal S_*$,  for any $y\in \mathbb B$, $y=\sum_{l=0}^s y_l$ with $\{y_l\in \mathbb B_l, l\in[\![0,s]\!]\}$ an $x\in \mathbb A$, $x=\sum_{n=0}^s x_n$ with $\{x_n\in \mathbb A_n, n\in[\![0,s]\!]\}$, such that
$ \mathcal T\left(x\right)=y$
can be defined component-wise by
$$
\left\{
\begin{array}{l}
x_0=\mathcal S_*( y_0), \\
\displaystyle 
\displaystyle x_l=\mathcal S_*\left( y_l-\sum_{k=0}^{l-1} \mathcal R_{lk}(x_k) \right),\ 
\forall l\in[\![1,s]\!].
\end{array}
\right.
$$
Hence summing the components for $l$ from $0$ to $s$ gives for $x=\sum_{l=0}^s x_l$
$$
x= \mathcal S_*\left(\sum_{l=0}^s y_l-\sum_{k=0}^{s-1}\sum_{l=k+1}^s \mathcal R_{lk}(x_k)\right)
\Rightarrow x= \mathcal S_* \left(y - \sum_{k=0}^{s-1}\mathcal R(x_k) \right)
\Rightarrow x = \mathcal S_* \left(y - \mathcal R\left(\sum_{k=0}^{s-1}x_k\right) \right).
$$
And since $\mathcal R(\mathbb A_s) = \{0_{\mathbb B}\}$, this actually gives $x = \mathcal S_*(y -\mathcal R(x))$, and, for any $k\in\mathbb N$, by $k$ iterated substitutions in the right hand side it gives
\begin{equation}\label{eq:it}
x = \mathcal S_*\left(\sum_{l=0}^k [-\mathcal R\circ\mathcal S_*]^l(y) \right)-\mathcal S_* [-\mathcal R\circ\mathcal S_*]^k\Big( \mathcal R(x)\Big).
\end{equation}
Besides, the operator $\mathcal R\circ\mathcal S_*: \mathbb B\to \mathbb B$ is nilpotent, and more precisely $ [-\mathcal R\circ\mathcal S_*]^{s+1}=0|_{\mathbb B\to \mathbb B}$; indeed, the range of $\mathcal R\circ\mathcal S_*$ restricted to any $ \bigoplus_{l=n}^{s} \mathbb B_l$ for any $n\in[\![0,s-1]\!]$ is a subset of $ \bigoplus_{l=n+1}^{s} \mathbb B_l$ and  the range of $\mathcal R\circ\mathcal S_*$ restricted to any $\mathbb B_s$ is $\{0_{\mathbb B}\}$.
This concludes the proof from \eqref{eq:it} with $k=s+1$.
\end{proof}
As a sanity check, it is also easy to verify directly that the proposed operator is indeed a right inverse of $\mathcal T$: for any $y\in \mathbb B$, then
$$
\mathcal T \left( \mathcal S_* \circ\sum_{j=0}^{s}[- \mathcal R \circ\mathcal S_*]^j (y)\right) 
=\mathcal T_* \left( \mathcal S_*  \left(\sum_{j=0}^{s}[- \mathcal R \circ\mathcal S_*]^j (y)\right) \right) 
+\mathcal R \left( \mathcal S_*  \left(\sum_{j=0}^{s}[- \mathcal R \circ\mathcal S_*]^j (y)\right) \right) ,
$$
$$
\Rightarrow \mathcal T \left( \mathcal S_* \circ\sum_{j=0}^{s}[- \mathcal R \circ\mathcal S_*]^j (y)\right) =\sum_{j=0}^{s}[- \mathcal R \circ\mathcal S_*]^j (y)
-\sum_{j=0}^{s}[- \mathcal R \circ\mathcal S_*]^{j+1} (y) = y-[- \mathcal R \circ\mathcal S_*]^{s+1} (y)  = y.
$$
\subsection{Quasi-Trefftz framework}\label{ssec:surj}
With the diagonal operator $\mathcal L_*$ playing a key role, the quasi-Trefftz operator $\mathcal D_p$ can be split  into the sum of its restriction to $\bigoplus_{\gamma}^p\widetilde{\mathbb P}_{\ell}$ and its restriction to $\bigoplus_0^{\gamma-1}\widetilde{\mathbb P}_{\ell}$. The former encompasses the essential properties of $\mathcal D_p$, as it satisfies the hypotheses of Theorem \ref{thm:surj}.
\begin{thm}\label{thm:qTsurj}
Given $(d,\gamma,p)\in\mathbb N^3$ such that $p>\gamma$ and $\mathbf x_0\in\mathbb R^d$,
consider
a linear differential operator of order at most equal to $\gamma$ with variable coefficients $\mathbf c = \{c_\mathbf j\in\mathcal C^{p-\gamma}(\mathbf x_0),\mathbf j\in\mathbb N_0^d, |\mathbf j| \leq \gamma\}$, $\displaystyle \mathcal L := \sum_{m=0}^\gamma \sum_{ |\mathbf j|=m} c_{\mathbf j} 
 \partial_{{\mathbf x}}^{\mathbf j}$, and the corresponding constant coefficient operator of order $\gamma$ defined by $\mathcal L_* = \sum_{ |\mathbf j|=\gamma} c_{\mathbf j}\left( {\mathbf x}_0  \right) \partial_{{\mathbf x}}^{\mathbf j}$.
\ Assume that the operator $\mathcal L_*|_{\mathbb P_p}:\mathbb P_p\to\mathbb P_{p-\gamma}$ is surjective. 
Then the quasi-Trefftz operator, $\mathcal D_p = T_{p-\gamma}\circ\mathcal L:\mathbb P_p\to\mathbb P_{p-\gamma}$, as well as $= \mathcal D_p|_{\bigoplus_{\gamma}^p\widetilde{\mathbb P}_{\ell}}$
 are both surjective.
 
 Moreover the operator $\mathcal S:\mathbb P_{p-\gamma}\to\mathbb P_{p}$ given by $\mathcal S = \mathcal S_* \circ\sum_{l=0}^{s}[- \mathcal R \circ\mathcal S_*]^l$ where $\mathcal R = \mathcal D_p-\mathcal L_*$ is a right inverse of $\mathcal D_p$.
 \end{thm}
 \begin{proof}
 The operator $\mathcal D_p $ is surjective if and only if for any $f\in\mathbb P_{p-\gamma}$ there exists $\Pi\in\mathbb P_p$ such that $\mathcal D_p \Pi = f$, 
 or equivalently if there exists $\{\Pi_\ell\in\widetilde{\mathbb P}_\ell, \ell\in[\![  0,p]\!]\}$ such that
 $$
\sum_{\ell = \gamma}^p \mathcal D_p ( \Pi_\ell)
=
f- \sum_{\ell = 0}^{\gamma-1}\mathcal D_p ( \Pi_\ell).
 $$
As announced earlier, the operator of interest is precisely the operator in the left-hand side of the previous equation, namely $\mathcal T = \mathcal D_p|_{\bigoplus_{\gamma}^p\widetilde{\mathbb P}_{\ell}}$.
 The operators $\mathcal L_*$ and the corresponding remainder $\mathcal R = \mathcal D_p-\mathcal L_*$ verify the following hypotheses.
\begin{enumerate}
\item
\begin{enumerate}
\item
For all $n\in[\![0,p-\gamma]\!]$, $\mathcal L_*( \widetilde{\mathbb P}_{n+\gamma})\subset  \widetilde{\mathbb P}_n$; 
indeed, $\mathcal L_*$ is a linear differential operator with constant coefficients and with only terms of order exactly $\gamma$.
\item The operator $\mathcal L_*$ is surjective by assumption.
\end{enumerate}
\item 
\begin{enumerate}
\item
For all $n\in[\![0,p-\gamma-1]\!]$, $\mathcal R(\widetilde{\mathbb P}_{n+\gamma})\subset \bigoplus_{l=n+1}^{p-\gamma} \widetilde{\mathbb P}_{n}$ from Lemma \ref{lmm:R}. 
\item
$\mathcal R(\widetilde{\mathbb P}_{p}) = \{0_{\mathbb P_{p-\gamma}}\}$;
indeed, for any $\Pi_p\in\widetilde{\mathbb P}_p$ then $\mathcal D_p(\Pi_p) = \mathcal L_*(\Pi_p)$, see Appendix \ref{app:newAdd}.
\end{enumerate}
\end{enumerate}
Considering
$s=p-\gamma$,
$\mathbb A_n = \widetilde{\mathbb P}_{n+\gamma}$,
and $\mathbb B_n = \widetilde{\mathbb P}_n$, as well as the decomposition $\mathcal D_p|_{\bigoplus_{\gamma}^p\widetilde{\mathbb P}_{\ell}} = \mathcal T_* + \mathcal R$ with $\mathcal T_* = \mathcal L_*$ and the corresponding remainder $\mathcal R = \mathcal D_p-\mathcal L_*$, these are precisely Theorem \ref{thm:surj}'s hypotheses, hence the theorem applies and the conclusion is that $\mathcal T = \mathcal D_p|_{\bigoplus_{\gamma}^p\widetilde{\mathbb P}_{\ell}}$ is surjective.

As a result, for any $f\in\mathbb P_{p-\gamma}$ and any $\{\Pi_\ell\in\widetilde{\mathbb P}_\ell, \ell\in[\![  0,\gamma-1]\!]\}$, there exists $\{\Pi_\ell\in\widetilde{\mathbb P}_\ell, \ell\in[\![  \gamma,p]\!]\}$ such that
 $$
\sum_{\ell = \gamma}^p \mathcal D_p ( \Pi_\ell)
=
f- \sum_{\ell = 0}^{\gamma-1}\mathcal D_p ( \Pi_\ell).
 $$. Following the initial remark of this proof, the operator $\mathcal D_p$ is then surjective.
\end{proof}

fully relies on the graded structure of polynomial spaces

$$
\mathcal D_p \left(\sum_{\ell = 0}^p \Pi_\ell\right)
=
\sum_{\ell = 0}^p\mathcal D_p ( \Pi_\ell)
=
\sum_{\ell = 0}^{\gamma-1}\mathcal D_p ( \Pi_\ell)
+
\sum_{\ell = \gamma}^p\underbrace{\mathcal D_p ( \Pi_\ell)}_{\in\cup_{k=\ell}^p \widetilde{\mathbb P}_k}
$$
\begin{lmm}\label{lmm:R}
Under the hypotheses of Theorem \ref{thm:qTsurj}, then 
for all $ N\in[\![0,p-\gamma-1]\!]$, $\mathcal R(\widetilde{\mathbb P}_{N+\gamma})\subset \bigoplus_{l=N-\gamma+1}^{p-\gamma} \widetilde{\mathbb P}_l$.
\end{lmm}
\begin{proof}
For any $N\in[\![\gamma,p-1]\!]$ and any $\Pi_N\in\widetilde{\mathbb P}_N$
$$
\mathcal R  (\Pi_N) = 
\sum_{m=\max(0,N-(p-\gamma))}^{\gamma-1} \sum_{ |\mathbf j|=m}  \sum_{\ell=N-m}^{p-\gamma} c_{\mathbf j,\ell-(N-m)} \partial_{\mathbf x}^{\mathbf j}\Pi_N
+\sum_{ |\mathbf j|=\gamma}  \sum_{\ell=N-\gamma+1}^{p-\gamma} c_{\mathbf j,\ell-(N-m)} \partial_{\mathbf x}^{\mathbf j}\Pi_N 
$$
with $c_{\mathbf j,\ell-(N-m)} \partial_{\mathbf x}^{\mathbf j}\Pi_N\in\widetilde{\mathbb P}_\ell$, see Appendix \ref{app:newexpl}.
Since $\ell\in[\![ N-\gamma+1 , p-\gamma ]\!]$ both in the triple sum and in the double sum, this concludes the proof.
\end{proof}

\section{Quasi-Trefftz space}\label{sec:ker}
Thanks to the surjectivity of the quasi-Trefftz operator, the quasi-Trefftz space, as the kernel of this operator, can be characterized as a direct sum of much simpler spaces.
This also provides about information the dimension of the quasi-Trefftz space.
\subsection{General framework}\label{ssec:GFker}
This section leverages the general framework of block-triangular linear operators to characterize the kernels of such operators.
\begin{thm}\label{thm:ker}
Consider $s\in\mathbb N$ as well as a linear operator $\mathcal T:\mathbb A\to \mathbb B$ between
 two graded vector spaces 
\begin{equation}
\mathbb A = \bigoplus_{n=0}^{s} \mathbb A_n\text{ and }\mathbb B = \bigoplus_{l=0}^{s} \mathbb B_l.
\end{equation}
Assuming that the operator can be decomposed as the sum of two operators $\mathcal T_*$ and $\mathcal R$ with
\begin{enumerate}
\item
\begin{enumerate}
\item\label{rangeT*bis}
for all $n\in[\![0,s]\!]$, $\mathcal T_*(\mathbb A_n)\subset \mathbb B_n$;
\item $\mathcal T_*:\mathbb A\to \mathbb B$ is surjective with a linear right inverse operator $\mathcal S_*$;
\item \label{rangeS}
for all $n\in[\![0,s]\!]$,  $\mathcal S_* (\mathbb B_n)\subset \mathbb A_n$;
\end{enumerate}
\item 
\begin{enumerate}
\item\label{rangeRbis}
for all $n\in[\![0,s-1]\!]$, $\mathcal R(\mathbb A_n)\subset \bigoplus_{l=n+1}^{s} \mathbb B_l$;
\item
$\mathcal R(\mathbb A_s) = \{0_{\mathbb B}\}$;
\end{enumerate}
\end{enumerate}
 then $\ker(\mathcal T)=[I_{\mathbb A}-\mathcal S\circ\mathcal R](\ker \mathcal T_*)$, where $I_{\mathbb A}$ denotes the identity operator defined on $\mathbb A$ and the operator $\mathcal S:= \mathcal S_* \circ\sum_{l=0}^{s}[- \mathcal R \circ\mathcal S_*]^l $ is the right inverse operator from Theorem \ref{thm:surj}.

Moreover, $\dim\ker(\mathcal T)=\dim\ker(\mathcal T_*)$.
\end{thm}
\begin{proof}
As a preliminary, a simple consequence of hypotheses \ref{rangeT*bis} and \ref{rangeRbis} reads:
\begin{equation}\label{eq:csq}
\forall n\in[\![ 0, s-1]\!],\mathcal T\left(\bigoplus_{l=n}^s \mathbb A_l\right)\subset \bigoplus_{l=n}^s \mathbb B_l.
\end{equation}
Moreover, for all $n\in[\![ 0, s-1]\!]$, for any $x\in\bigoplus_{l=n}^s \mathbb A_l$ there exists a unique $(y_n,z_n)\in \mathbb A_n\times\bigoplus_{l=n+1}^s \mathbb A_l$ such that $x = y_n+z_n$, and a further consequence of these hypotheses and \eqref{eq:csq} reads
$$
\mathcal T (x)  = \mathcal T_* (y_n) + \mathcal R(y_n)+\mathcal T(z_n)
\text{ with }
\left\{\begin{array}{l}
 \mathcal T_* (y_n) \in \mathbb B_n ,\\
  \mathcal R(y_n)+\mathcal T(z_n) \in \bigoplus_{l=n+1}^s \mathbb B_l.
\end{array}\right.
$$
In particular, $\mathcal T(x) = 0_{\bigoplus_{l=n}^s \mathbb B_l}$ implies both that $\mathcal T_* y_n = 0_{\mathbb B_n}$ and that $\mathcal T(z_n) = - \mathcal R(y_n)$. In other words, for any $x\in\ker \mathcal T|_{\bigoplus_{l=n}^s \mathbb A_l}$, there exists $(y_n,w_n)\in \ker \mathcal T_*|_{\mathbb A_n}\times \ker\mathcal T|_{\bigoplus_{l=n+1}^s \mathbb A_l}$ such that $x = y_n-\mathcal S\mathcal R (y_n) +w_n$.
Hence this relates the kernel of $\mathcal T$ restricted to $\bigoplus_{l=n}^s \mathbb A_l$ to the kernel of $\mathcal T$ restricted to $\bigoplus_{l=n+1}^s \mathbb A_l$.

Combining these statements for increasing values of $n\in[\![ 0, s-1]\!]$ then shows that for any $x\in\ker \mathcal T$, there exists $w_{n-1}\ker\mathcal T|_{\mathbb A_s}$ and $\{y_n\in \ker \mathcal T_*|_{\mathbb A_n}, n\in[\![ 0, s-1]\!]\}$ such that $x = y_n-\mathcal S\mathcal R (y_n) +w_n$
$$
\forall x\in\ker \mathcal T, 
\exists w_{s-1}\in\ker\mathcal T|_{\mathbb A_s}\text{ and }\{y_n\in \ker \mathcal T_*|_{\mathbb A_n}, n\in[\![ 0, s-1]\!]\}
\text{ such that }x = \sum_{n=0}^{s-1} [I_{\mathbb A}-\mathcal S\circ\mathcal R](y_n) + w_{s-1} .
$$
Moreover, since $\ker\mathcal T|_{\mathbb A_s} = \ker\mathcal T_*|_{\mathbb A_s}$ and $\mathcal R(\mathbb A_s) = \{0_{\mathbb B}\}$, then $w_{s-1}\in \ker\mathcal T_*|_{\mathbb A_s}$ and $w_{s-1} =  [I_{\mathbb A}-\mathcal S\circ\mathcal R](w_{s-1})$. This  gives
$$
\forall x\in\ker \mathcal T, 
\exists\{y_n\in \ker \mathcal T_*|_{\mathbb A_n}, n\in[\![ 0, s]\!]\}
\text{ such that }x = \sum_{n=0}^{s-1} [I_{\mathbb A}-\mathcal S\circ\mathcal R](y_n),
$$
$$
\Leftrightarrow
\forall x\in\ker \mathcal T, 
y\in\ker \mathcal T_*
\text{ such that }x =  [I_{\mathbb A}-\mathcal S\circ\mathcal R](y)  \text{ as } \ker \mathcal T_* = \bigoplus_{n=0}^s \ker \mathcal T_*|_{\mathbb A_n},
$$
as the last identity follows from hypothesis \ref{rangeT*bis}. 
Hence $\ker(\mathcal T)\subset[I_{\mathbb A}-\mathcal S\circ\mathcal R](\ker \mathcal T_*)$.

Conversely, for any $y\in[I_{\mathbb A}-\mathcal S\circ\mathcal R](\ker \mathcal T_*)$, there exists $x\in \ker \mathcal T_*$ such that $y = [I_{\mathbb A}-\mathcal S\circ\mathcal R]x$, so that $\mathcal Ty = \mathcal T(x-\mathcal S(\mathcal R(x))) =0 $ since $\mathcal S$ is a right inverse of $\mathcal T$, $\mathcal T = \mathcal T_*+\mathcal R$ and $\mathcal T_*(x) = 0 $.
Hence $[I_{\mathbb A}-\mathcal S\circ\mathcal R](\ker \mathcal T_*)\subset\ker(\mathcal T)$, which in turn shows the desired identity: $\ker(\mathcal T)=[I_{\mathbb A}-\mathcal S\circ\mathcal R](\ker \mathcal T_*)$.

Finally, in order to prove that the dimensions of the kernels of $\mathcal T$ and $\mathcal T_*$ are equal, thanks to the rank-nullity theorem for the operator $[I_{\mathbb A}-\mathcal S\circ\mathcal R]|_{\ker \mathcal T_*}$ it is sufficient to prove that the kernel of this operator is trivial.
For any $y\in\ker [I_{\mathbb A}-\mathcal S\circ\mathcal R]|_{\ker \mathcal T_*}$ then $y=\mathcal S(\mathcal R(y))$. Assuming that $y\neq 0_{\mathbb A}$, then by Lemma \ref{lmm:SR}:
\begin{itemize}
\item either $y\in \mathbb A_s$, so $\mathcal S(\mathcal R (y)) = 0_{\mathbb A}$;
\item or there exists $N\in[\![0,s-1]\!]$ such that $y\in\bigoplus_{n=N}^{s} \mathbb A_n$ and $y\notin\bigoplus_{n=N+1}^{s} \mathbb A_n$, so $\mathcal S(\mathcal R(y))\in  \bigoplus_{n=N+1}^{s} \mathbb A_n$;
\end{itemize}
in either case this is a contradiction since $y=\mathcal S(\mathcal R(y))\neq 0_{\mathbb A}$. Hence any $y\in\ker [I_{\mathbb A}-\mathcal S\circ\mathcal R]|_{\ker \mathcal T_*}$ is zero.

This concludes the proof.
\end{proof}

\begin{lmm}\label{lmm:SR}
Under the hypotheses of Theorem \ref{thm:ker}, then 
$$\mathcal S \circ \mathcal R(\mathbb A_s)\subset \{0_{\mathbb A}\} 
\text{ while }
\forall N\in[\![0,s-1]\!],\mathcal S \circ \mathcal R\left(\bigoplus_{n=N}^{s} \mathbb A_n\right)\subset \bigoplus_{n=N+1}^{s} \mathbb A_n.$$
\end{lmm}
\begin{proof}
First, for any $N\in[\![0,s-1]\!]$  the image of ${ \bigoplus_{n=N}^{s} \mathbb A_n}$ by $\mathcal R$ is $ \bigoplus_{n=N+1}^{s} \mathbb B_n$. 
Next, 
 $$
 \left\{\begin{array}{l}
[\mathcal R\circ\mathcal S_*](\mathbb B_s)=\{0_{\mathbb B}\},
\\ 
\forall N\in[\![0,s-1]\!],
\\
\displaystyle{} [\mathcal R\circ\mathcal S_*]\left( \bigoplus_{n=N}^{s} \mathbb B_n\right)\subset\bigoplus_{n=N+1}^{s} \mathbb B_n,
 \end{array}\right.
 \Rightarrow \forall l\in[\![1,s]\!],
 \left\{\begin{array}{l}
\displaystyle
\forall N\in[\![s-l+1,s]\!], [\mathcal R\circ\mathcal S_*]^l \left( \bigoplus_{n=N}^{s} \mathbb B_n\right)=\{0_B\},\\
\displaystyle
\forall N\in[\![0,s-l]\!],[\mathcal R\circ\mathcal S_*]^l\left( \bigoplus_{n=N}^{s} \mathbb B_n\right)\subset\bigoplus_{k=N+l}^{s} \mathbb B_k.
 \end{array}\right.
 $$
 This implies
 $$
 \forall N\in[\![0,s-2]\!], \sum_{l=1}^{s}[- \mathcal R \circ\mathcal S_*]^l \left({\bigoplus_{n=N+1}^{s} \mathbb B_n}\right)\subset \bigoplus_{n=N+2}^{s} \mathbb B_n
 $$
 and in turn
 $$
\forall N\in[\![0,s-1]\!], \sum_{l=0}^{s}[- \mathcal R \circ\mathcal S_*]^l \left( {\bigoplus_{n=N+1}^{s} \mathbb B_n}\right)\subset \bigoplus_{n=N+1}^{s} \mathbb B_n.
 $$
 Combined with the hypothesis on $\mathcal S_*$, this gives that
$$
\forall N\in[\![0,s-1]\!], \mathcal S \circ \mathcal R\left({\bigoplus_{n=N}^{s} \mathbb A_n}\right)\subset\bigoplus_{k=N+1}^{s} \mathbb A_k.
$$
As for the case $N=s$, it directly follows from the facts that $\mathcal R(\mathbb A_s) = \{0_{\mathbb B}\}$ and $\mathcal S_*(0_{\mathbb B}) = 0_{\mathbb A}$.
\end{proof}

\subsection{Quasi-Trefftz framework}\label{ssec:ker} 
Here again, the operator $\mathcal L_*$ plays a key role as the diagonal part of $\mathcal D_p$ restricted to $\bigoplus_{\gamma}^p\widetilde{\mathbb P}_{\ell}$, and the latter satisfies the hypotheses of Theorem \ref{thm:ker}.
\begin{thm}\label{thm:qTker}
Given $(d,\gamma,p)\in\mathbb N^3$ such that $p>\gamma$ and $\mathbf x_0\in\mathbb R^d$,
consider
a linear differential operator of order at most equal to $\gamma$ with variable coefficients $\mathbf c = \{c_\mathbf j\in\mathcal C^{p-\gamma}(\mathbf x_0),\mathbf j\in\mathbb N_0^d, |\mathbf j| \leq \gamma\}$, $\displaystyle \mathcal L := \sum_{m=0}^\gamma \sum_{ |\mathbf j|=m} c_{\mathbf j} 
 \partial_{{\mathbf x}}^{\mathbf j}$, and the corresponding constant coefficient operator of order $\gamma$ defined by $\mathcal L_* = \sum_{ |\mathbf j|=\gamma} c_{\mathbf j}\left( {\mathbf x}_0  \right) \partial_{{\mathbf x}}^{\mathbf j}$, so that for all $n\in\mathbb N_0$, $\mathcal L_* :\widetilde{\mathbb P}_{n+\gamma}\to\widetilde{\mathbb P}_n$.
 Assume that the operator $\mathcal L_*|_{\mathbb P_p}:\mathbb P_p\to\mathbb P_{p-\gamma}$ is surjective with a linear right inverse operator $\mathcal S_*:\mathbb P_{p-\gamma}\to \mathbb P_{p}$ satisfying for all $n\in[\![0,p-\gamma]\!]$,  $\mathcal S_* (\widetilde{\mathbb P}_n)\subset \widetilde{\mathbb P}_{n+\gamma}$. 
For the corresponding quasi-Trefftz operator, $\mathcal D_p = T_{p-\gamma}\circ\mathcal L:\mathbb P_p\to\mathbb P_{p-\gamma}$,
define $\mathcal T = \mathcal D_p|_{\bigoplus_{\gamma}^p\widetilde{\mathbb P}_{\ell}}$.
Then $$\ker(\mathcal T)=[I_{\mathbb P_p}-\mathcal S\circ\mathcal R](\ker \mathcal L_*),$$ where $I_{\mathbb P_p}$ denotes the identity operator defined on $\mathbb P_p$ and the operator $\mathcal S:= \mathcal S_* \circ\sum_{l=0}^{s}[- \mathcal R \circ\mathcal S_*]^l $ is the right inverse operator from Theorem \ref{thm:qTsurj}.

Moreover, $\dim\ker(\mathcal D_p)=\dim\ker(\mathcal L_*)+ \dim \bigoplus_0^{\gamma-1}\widetilde{\mathbb P}_{\ell}$.
 \end{thm}
 \begin{proof}
 As a preliminary remark, the set of hypotheses of Theorem \ref{thm:qTker} include all hypotheses of Theorem \ref{thm:qTsurj}, hence both Theorem \ref{thm:qTsurj} and Lemma \ref{lmm:R} apply.
 
Here again, the operator of interest is $\mathcal T = \mathcal D_p|_{\bigoplus_{\gamma}^p\widetilde{\mathbb P}_{\ell}}$.
 The operators $\mathcal L_*$ and $\mathcal R = \mathcal D_p-\mathcal L_*$, the corresponding remainder,  verify the following hypotheses.
\begin{enumerate}
\item
\begin{enumerate}
\item
for all $n\in[\![0,s]\!]$, $\mathcal L_*(\widetilde{\mathbb P}_{n+\gamma})\subset \widetilde{\mathbb P}_n$;
\item $\mathcal L_*:\bigoplus_{\gamma}^p\widetilde{\mathbb P}_{\ell}\to \bigoplus_{0}^{p-\gamma}\widetilde{\mathbb P}_{\ell}$ is surjective by hypothesis;
\item 
for all $n\in[\![0,s]\!]$ $\mathcal S_*|_{\widetilde{\mathbb P}_n}:\widetilde{\mathbb P}_n\to \widetilde{\mathbb P}_{n+\gamma}$ by hypothesis;
\end{enumerate}
\item 
\begin{enumerate}
\item
for all $n\in[\![0,p-\gamma-1]\!]$, $\mathcal R(\widetilde{\mathbb P}_{n+\gamma})\subset \bigoplus_{l=n+1}^{p-\gamma} \widetilde{\mathbb P}_l$ from Lemma \ref{lmm:R};
\item
$\mathcal R(\widetilde{\mathbb P}_{p+\gamma}) = \{0_{\mathbb P}\}$;
\end{enumerate}
\end{enumerate}
Considering
$s=p-\gamma$,
for all $n\in[\![0,p-\gamma]\!]$ $\mathbb A_n = \widetilde{\mathbb P}_{n+\gamma}$ and $\mathbb B_n = \widetilde{\mathbb P}_n$, 
as well as the decomposition $\mathcal D_p|_{\bigoplus_{\gamma}^p\widetilde{\mathbb P}_{\ell}} = \mathcal T_* + \mathcal R|_{\bigoplus_{\gamma}^p\widetilde{\mathbb P}_{\ell}}$ with $\mathcal T_* = \mathcal L_*$ and the corresponding remainder $\mathcal R = \mathcal D_p-\mathcal L_*$, these are precisely Theorem \ref{thm:ker}'s hypotheses, hence the theorem applies. The conclusion is that for $\mathcal T := \mathcal D_p|_{\bigoplus_{\gamma}^p\widetilde{\mathbb P}_{\ell}}$  then $\ker(\mathcal T)=[I_{\bigoplus_{\gamma}^p\widetilde{\mathbb P}_{\ell}}-\mathcal S\circ\mathcal R](\ker \mathcal L_*)$ and $\dim\ker(\mathcal T)=\dim\ker(\mathcal L_*)$.

Finally, the rank nullity theorem leads to
$$
\dim\bigoplus_{\ell=0}^p\widetilde{\mathbb P}_{\ell} = \dim \ker(\mathcal D_p )+rk(\mathcal D_p )
\text{ and }
\dim \bigoplus_{\ell=\gamma}^p\widetilde{\mathbb P}_{\ell} = 
\dim\ker(\mathcal L_*) 
+
 rk\left(\mathcal D_p|_{\bigoplus_{\ell=\gamma}^p\widetilde{\mathbb P}_{\ell}} \right).
$$
Moreover, from Theorem \ref{thm:qTsurj} both $\mathcal D_p$ and $\mathcal T := \mathcal D|_{\bigoplus_{\ell=\gamma}^p\widetilde{\mathbb P}_{\ell}}$ are surjective: 
$$
 rk\left(\mathcal D_p|_{\bigoplus_{\ell=\gamma}^p\widetilde{\mathbb P}_{\ell}} \right) =\dim \mathbb P_{p-\gamma}
 \text{ and }
 rk\left(\mathcal D_p \right) =\dim \mathbb P_{p-\gamma}
$$
 This then concludes the proof as it implies
 $$
 \dim\bigoplus_{\ell=0}^p\widetilde{\mathbb P}_{\ell} - \dim \ker(\mathcal D_p )
 =
 \dim \bigoplus_{\ell=\gamma}^p\widetilde{\mathbb P}_{\ell} -\dim\ker(\mathcal L_*) .
 $$
\end{proof}

\section{Quasi-Trefftz basis}\label{sec:bas}
Thanks to the surjectivity of the quasi-Trefftz operator and the characterization of its kernel, a quasi-Trefftz basis can be constructed thanks to a general procedure.

\subsection{General framework}\label{ssec:GFbas}
The section describes the set of solutions of block-triangular systems thanks to a forward substitution procedure, 
as well as the characterization of the set of solutions.
The procedure leverages a sequence of linear systems. 
\begin{prop}\label{prop:kers}
 Given a set of indices $I$, two vector spaces $\mathbb F$ and $\mathbb G$ with $\dim \mathbb F>\dim \mathbb G$, a linear operator $\widetilde{\mathcal T}:\mathbb F\to \mathbb G$, assume that $\mathbb F=\mathbb U\oplus \mathbb V$,  $\widetilde{\mathcal T}$ is surjective, and $\widetilde{\mathcal T}|_{\mathbb V}$ is invertible.
Then 
 given any $Y\in \mathbb B$, the set of solutions to $\widetilde{\mathcal T} X=Y$ is $\{U+[\widetilde{\mathcal T}|_{\mathbb V}]^{-1}(Y-\widetilde{\mathcal T}(U)),\forall U\in \mathbb U \}$.
\end{prop}
\begin{proof}
First, $\{X_i-[\widetilde{\mathcal T}|_{\mathbb V}]^{-1}(\widetilde{\mathcal T}(X_i)),i\in I\}$ is clearly a subset of $\ker(\widetilde{\mathcal T})$. 
Second, it is linearly independent according to Lemma \ref{lem:linind}.
Moreover, $\dim\ker\widetilde{\mathcal T}=\dim \mathbb U$, since $\dim \mathbb F = \dim\ker\widetilde{\mathcal T}+rk(\widetilde{\mathcal T})$, $rk\widetilde{\mathcal T}=rk\widetilde{\mathcal T}|_{\mathbb V}$ and $\dim \mathbb F = \dim \mathbb U+\dim \mathbb V$.
Hence, $\{X_i-[\widetilde{\mathcal T}|_{\mathbb V}]^{-1}(\widetilde{\mathcal T}(X_i))\}\subset\ker\widetilde{\mathcal T}$ is a set of $\dim\ker\widetilde{\mathcal T}$ linearly independent vectors, so it forms a basis of $\ker\widetilde{\mathcal T}$.
Finally, the candidate set is the sum of a particular solution $[\widetilde{\mathcal T}|_{\mathbb V}]^{-1}(Y)$ and the kernel of $\widetilde{\mathcal T}$. This concludes the proof.
\end{proof}
\begin{lmm}\label{lem:linind}
Given a set of indices $I$, two vector spaces $\mathbb U$ and $\mathbb V$ as well as $\{(X_i,Y_i),i\in I\}\subset \mathbb U\times \mathbb V$, 
if $\{X_i,i\in I\}\subset \mathbb U$ is linearly independent, then $\{X_i+Y_i,i\in I\}\subset \mathbb U\oplus \mathbb V$ is also linearly independent.
\end{lmm}
\begin{proof}
For any $\{\lambda_i\in\mathbb K,i\in I\}$, then
$$
\left[\sum_{i\in I}\lambda_i(X_i+Y_i)=0\right]
\Leftrightarrow
\left[\left(\sum_{i\in I}\lambda_iX_i\right)+\left(\sum_{i\in I}\lambda_iY_i\right)=0\right]
\Rightarrow 
\left[\sum_{i\in I}\lambda_iX_i = 0\right]
$$
since $\mathbb U\oplus \mathbb V$, and this implies that all $\lambda_i$s are zero as the set $\{X_i,i\in I\}$ is linearly independent. This concludes the proof.
\end{proof}
Thanks to Proposition \ref{prop:kers}, it is now possible to construct a basis for the set of solutions. 
\begin{thm}\label{thm:basis}
Consider $s\in\mathbb N$ as well as a linear operator $\mathcal T:\mathbb A\to \mathbb B$ between
 two graded vector spaces 
\begin{equation}
\mathbb A = \bigoplus_{n=0}^{s} \mathbb A_n\text{ and }\mathbb B = \bigoplus_{l=0}^{s} \mathbb B_l.
\end{equation}
Assuming that the operator can be decomposed as the sum of two operators $\mathcal T_*$ and $\mathcal R$ with
\begin{enumerate}
\item
\begin{enumerate}
\item
for all $n\in[\![0,s]\!]$, $\mathcal T_*(\mathbb A_n)\subset \mathbb B_n$;
\item for all $n\in[\![0,s]\!]$, $\mathbb A_n=\mathbb U_n\oplus \mathbb V_n$ and  $\mathcal T_*|_{\mathbb V_n}$ is bijective, 
its inverse denoted $\mathcal S_n$ for convenience;
\end{enumerate}
\item 
\begin{enumerate}
\item
for all $n\in[\![0,s-1]\!]$, $\mathcal R(\mathbb A_n)\subset \bigoplus_{l=n+1}^{s} \mathbb B_l$;
 for convenience $\mathcal R|_{\mathbb A_n}:\mathbb A_n\to \bigoplus_{l=n+1}^{s} \mathbb B_l$ is decomposed as the sum of $s-n$ operators $\mathcal R_{ln}:\mathbb A_n\to \mathbb B_l$ for all $l\in[\![n+1,s]\!]$:
$$
\forall x_n\in \mathbb A_n, \mathcal R(x_n) = \sum_{l=n+1}^s \mathcal R_{ln}(x_n);
$$
\item
$\mathcal R(\mathbb A_s) = \{0_{\mathbb B}\}$;
\end{enumerate}
\end{enumerate}
then, given any $b\in\mathbb B$, Algorithm \ref{algo:ker} computes the components $\{x_n\in \mathbb A_n,n\in[\![0,s]\!]\}$ of a generic element $X\in\mathbb A$ is the set of solutions to $\mathcal T(X)=b$:

\begin{algorithm}[!ht]
\label{algo:ker}
\DontPrintSemicolon
	{\bf Input} 
	$\displaystyle b = \sum_{k=0}^sb_k$ with  $b_k\in \mathbb B_k\ \forall k \in[\![0,s]\!]$
	and
	$\displaystyle U = \sum_{k=0}^su_k$ with $u_k\in \mathbb U_k\ \forall k \in[\![0,s]\!]$

    $x_0 \gets u_0+\mathcal S_0(b_0-\mathcal T_*(u_0))$
    
    \For{$k \in[\![1,s]\!]$\,}{
    	$\displaystyle x_k \gets  u_k+\mathcal S_k\left(b_k- \sum_{j=0}^{k-1} \mathcal R_{kj}(x_j) - \mathcal T_*(u_k)\right)$
    }
    
    {\bf Output} $\displaystyle X = \sum_{k=0}^s x_k$
    
\caption{Construction of a generic element of the kernel}
\end{algorithm}
\FloatBarrier

Moreover, given $b\in{\mathbb B}$, for any basis $\{Y_i,i\in I\}$ of $\mathbb U = \oplus_{k=0}^s \mathbb U_k$ for some set of indices $I$, consider the corresponding  $\{X_i,i\in I\}$ where, for each $i\in I$, $X_i$ is the output of Algorithm \ref{algo:ker} for the input $(b,Y_i)$. Then $\{X_i+Y_i,i\in I\}$ forms a basis of the set of solutions. 
\end{thm}
\begin{proof}

Algorithm \ref{algo:ker}  simply leverages Proposition \ref{prop:kers} at each step in the forward substitution procedure with the operator $\widetilde{\mathcal T}=\mathcal T_*$ from $\mathbb F=\mathbb A_k$, with $\mathbb A_k=\mathbb U_k\oplus\mathbb V_k$, to $\mathbb G=\mathbb B_k$. 

For convenience, denote $\mathbb U = \oplus_{k=0}^s \mathbb U_k$ and $\mathbb V = \oplus_{k=0}^s \mathbb V_k$.
Then, given $b=0_\mathbb B$, for any basis $\{Y_i,i\in I\}$ of $\mathbb U$, consider the corresponding  $\{X_i,i\in I\}$ where, for each $i\in I$, $X_i$ is the output of Algorithm \ref{algo:ker} for the input $Y_i$. Then $\{X_i+Y_i,i\in I\}\subset\ker\mathcal T$
 by construction and this set is linearly independent according to Lemma \ref{lem:linind}, so it spans a space of dimension $\dim\mathbb U$. Moreover, since the set of hypotheses of Theorem \ref{thm:basis} either include or imply all those of Theorem \ref{thm:ker}, Theorem \ref{thm:ker} applies and in particular $\dim\ker(\mathcal T)=\dim\ker(\mathcal T_*)$. 
On the other hand:
$$
\begin{array}{rll}
 \dim\ker(\mathcal T_*) 
 &\displaystyle= \sum_{n=0}^s \dim\ker(\mathcal T_*|_{\mathbb A_n}), 
 &\displaystyle \text{ as }\mathbb A = \bigoplus_{n=0}^{s} \mathbb A_n,\\
 &\displaystyle= \sum_{n=0}^s\dim\mathbb A_n -rk \mathcal T_*|_{\mathbb A_n} ,
 &\text{ from the rank nullity theorem for }\mathcal T_*|_{\mathbb A_n},\\
 &\displaystyle\leq \sum_{n=0}^s\dim\mathbb A_n -rk \mathcal T_*|_{\mathbb V_n} ,
 &\text{ since }\mathbb V_n\subset{\mathbb A_n} \text{ so }rk\mathcal T_*|_{\mathbb A_n}\geq  rk\mathcal T_*|_{\mathbb V_n},\\ 
 &\displaystyle\leq \sum_{n=0}^s\dim\mathbb A_n -\dim{\mathbb V_n} ,
 &\text{ since } \mathcal T_*|_{\mathbb V_n} \text{ is bijective},\\
 &\displaystyle\leq \sum_{n=0}^s\dim\mathbb U_n ,
 &\text{ since }\mathbb A_n=\mathbb U_n\oplus \mathbb V_n,\\
 &\displaystyle\leq \dim\mathbb U ,
 &\text{ by definition of }\mathbb U,\\
\end{array}
$$
and therefore $ \dim\ker(\mathcal T_*)  = \dim\mathbb U$.
As a result, $\dim\mathbb U = \dim\ker\mathcal T$; and therefore $\{X_i+Y_i,i\in I\}$ is s a basis for $\ker\mathcal T$.

Finally, given $b\in\mathbb B$, for any basis $\{Y_i,i\in I\}$ of $\mathbb U$, consider the corresponding  $\{X_i,i\in I\}$ where, for each $i\in I$, $X_i$ is the output of Algorithm \ref{algo:ker} for the input $Y_i$. Then  $ \sum_{k=0}^s \mathcal S_k b_k$ is a particular solution, and $\{X_i+Y_i,i\in I\}$ forms a basis of the set of solutions to $\mathcal T(x)=b$.
\end{proof}

\subsection{Quasi-Trefftz framework}\label{ssec:bas}
As $\mathcal D_p$ restricted to $\bigoplus_{\gamma}^p\widetilde{\mathbb P}_{\ell}$ satisfies the hypotheses of Theorem \ref{thm:basis}, the constructions of a quasi-Trefftz basis follows the forward substitution procedure, and each iteration $k$ requires to solve a system on a space of homogeneous polynomials $\widetilde{\mathbb P}_k$ for the operator $\mathcal L_*$.
\begin{thm}\label{thm:qTbasis}
Given $(d,\gamma,p)\in\mathbb N^3$ such that $p>\gamma$ and $\mathbf x_0\in\mathbb R^d$,
consider
a linear differential operator of order at most equal to $\gamma$ with variable coefficients $\mathbf c = \{c_\mathbf j\in\mathcal C^{p-\gamma}(\mathbf x_0),\mathbf j\in\mathbb N_0^d, |\mathbf j| \leq \gamma\}$, $\displaystyle \mathcal L := \sum_{m=0}^\gamma \sum_{ |\mathbf j|=m} c_{\mathbf j} 
 \partial_{{\mathbf x}}^{\mathbf j}$, and the corresponding constant coefficient operator of order $\gamma$ defined by $\mathcal L_* = \sum_{ |\mathbf j|=\gamma} c_{\mathbf j}\left( {\mathbf x}_0  \right) \partial_{{\mathbf x}}^{\mathbf j}$, so that for all $n\in\mathbb N_0$, $\mathcal L_* |_{\widetilde{\mathbb P}_{n+\gamma}}:\widetilde{\mathbb P}_{n+\gamma}\to\widetilde{\mathbb P}_n$.
 Assume that 
for all $n\in[\![\gamma,p]\!]$, $\widetilde{\mathbb P}_n=\mathbb U_n\oplus \mathbb V_n$ and  $\mathcal L_*|_{\mathbb V_n}:\mathbb V_n\to\widetilde{\mathbb P}_{n-\gamma}$ is bijective, 
its inverse denoted $\mathcal S_n$. 
If $\mathcal R:=\mathcal L-\mathcal L_*$, for convenience 
for all $n\in[\![0,p-\gamma]\!]$ $\mathcal R|_{\widetilde{\mathbb P}_{n+\gamma}}:\widetilde{\mathbb P}_{n+\gamma}\to \bigoplus_{l=n+1}^{p-\gamma} \widetilde{\mathbb P}_l$ is decomposed as the sum of $p-\gamma-n$ operators $\mathcal R_{l,{n+\gamma}}:\widetilde{\mathbb P}_{n+\gamma}\to \widetilde{\mathbb P}_l$ for all $ l\in[\![n+1,p-\gamma]\!]$,
$$
\forall n\in[\![0,p-\gamma]\!],
\forall x_{n+\gamma}\in \widetilde{\mathbb P}_{n+\gamma}, \mathcal R(x_{n+\gamma}) = \sum_{l=n+1}^{p-\gamma} \mathcal R_{l,n+\gamma}(x_{n+\gamma})
,
$$ and
for all $n\in[\![0,\gamma-1]\!]$ $\mathcal R|_{\widetilde{\mathbb P}_{n}}:\widetilde{\mathbb P}_{n}\to \bigoplus_{l=0}^{p-\gamma} \widetilde{\mathbb P}_l$ is decomposed as the sum of $p-\gamma$ operators $\mathcal R_{l,{n}}:\widetilde{\mathbb P}_{n}\to \widetilde{\mathbb P}_l$ for all $ l\in[\![0,p-\gamma]\!]$:
$$
\forall n\in[\![0,\gamma-1]\!],
\forall x_{n}\in \widetilde{\mathbb P}_{n}, \mathcal R(x_{n}) = \sum_{l=0}^{p-\gamma} \mathcal R_{l,n}(x_{n}).
$$ 
Then the Algorithm \ref{algo:qTker} computes the components $\{x_n\in \widetilde{\mathbb P}_n,n\in[\![0,p]\!]\}$ of a generic element $X\in\mathbb Q\mathbb T_p$:

\begin{algorithm}[!ht]
\label{algo:qTker}
\DontPrintSemicolon
	{\bf Input} 
	$\displaystyle U = \sum_{k=0}^{p}u_k$ with $u_{k}\in \mathbb U_k\ \forall k \in[\![\gamma,p]\!]$ and $u_k\in \widetilde{\mathbb P}_k \ \forall k \in[\![0,\gamma-1]\!]$
	
    \For{$k \in[\![0,\gamma-1]\!]$\,}{
    $x_k \gets u_k$
    }
        
    \For{$k \in[\![\gamma,p]\!]$\,}{
    	$\displaystyle x_{k} \gets u_{k}-\mathcal S_k\left( \sum_{j=0}^{k-1} \mathcal R_{k,j}(x_j) + \mathcal L_*(u_{k})\right)$
    }
    
    {\bf Output} $\displaystyle X = \sum_{k=0}^s x_k$
    
\caption{Construction of a generic element of the kernel}
\end{algorithm}
\FloatBarrier

Moreover, for any basis $\{Y_i,i\in I\}$ of $\mathbb U =\left(  \bigoplus_{n=0}^{\gamma-1} \widetilde{\mathbb P}_n\right) \oplus\left(\bigoplus_{k=\gamma}^{p} \mathbb U_k\right)$ for some set of indices $I$, consider the corresponding  $\{X_i,i\in I\}$ where, for each $i\in I$, $X_i$ is the output of Algorithm \ref{algo:qTker} for the input $Y_i$. Then $\{X_i+Y_i,i\in I\}$ forms a basis of $\mathbb Q\mathbb  T_p$.  And as a consequence, $\dim\mathbb Q\mathbb  T_p = \dim \mathbb P_{p}- \dim \mathbb P_{p-\gamma}$.
 \end{thm}

\begin{proof}
As a preliminary remark,
for any $\Pi\in\mathbb P_p$ with $\Pi=\Pi_1+\Pi_2$, $(\Pi_1,\Pi_2)\in \mathbb P_{\gamma-1}\times \bigoplus_\gamma^{p}\widetilde{\mathbb P}_{\ell} $, then
$$
\mathcal D_p \Pi = 0_{\mathbb P} \Leftrightarrow  \mathcal D_p|_{\bigoplus_{\gamma}^p\widetilde{\mathbb P}_{\ell}}
 \Pi_2 = -\mathcal D_p \Pi_1.
$$
Moreover, the set of hypotheses of Theorem \ref{thm:qTbasis} include all hypotheses of Theorem \ref{thm:qTker}, hence Theorem \ref{thm:qTker} applies, in particular
$\dim\ker(\mathcal D_p)=\dim\ker(\mathcal L_*)+ \dim \bigoplus_0^{\gamma-1}\widetilde{\mathbb P}_{\ell}$.
 
Given the operator $\mathcal T = \mathcal D_p|_{\bigoplus_{\gamma}^p\widetilde{\mathbb P}_{\ell}}$.
the operators of interest are $\mathcal L_*$ and $\mathcal R = \mathcal D_p-\mathcal L_*$, the corresponding remainder,  verify the following hypotheses:
\begin{enumerate}
\item
\begin{enumerate}
\item
for all $n\in[\![0,p-\gamma]\!]$, $\mathcal L_*(\widetilde{\mathbb P}_{n+\gamma})\subset \widetilde{\mathbb P}_n$;
\item for all $n\in[\![0,p-\gamma]\!]$,
$\widetilde{\mathbb P}_n=\mathbb U_n\oplus \mathbb V_n$ and $\mathcal L_*|_{\mathbb V_n}$ is bijective by hypothesis, its inverse denoted $\mathcal S_n$;
\end{enumerate}
\item 
\begin{enumerate}
\item
for all $n\in[\![0,p-\gamma-1]\!]$, $\mathcal R(\widetilde{\mathbb P}_{n+\gamma})\subset \bigoplus_{l=n+1}^{p-\gamma} \widetilde{\mathbb P}_l$ from Lemma \ref{lmm:R}, with $\mathcal R|_{\widetilde{\mathbb P}_{n+\gamma}}$ decomposed into corresponding operators $\mathcal R_{l,{n+\gamma}}$;
\item
$\mathcal R(\widetilde{\mathbb P}_{p+\gamma}) = \{0_{\mathbb P}\}$;
\end{enumerate}
\end{enumerate}
Considering
$s=p-\gamma$,
for all $n\in[\![0,p-\gamma]\!]$ $\mathbb A_n = \widetilde{\mathbb P}_{n+\gamma}$ and $\mathbb B_n = \widetilde{\mathbb P}_n$, 
as well as the decomposition $\mathcal D_p|_{\bigoplus_{\gamma}^p\widetilde{\mathbb P}_{\ell}} = \mathcal T_* + \mathcal R|_{\bigoplus_{\gamma}^p\widetilde{\mathbb P}_{\ell}}$ with $\mathcal T_* = \mathcal L_*$ and the corresponding remainder $\mathcal R = \mathcal D_p-\mathcal L_*$, these are precisely Theorem \ref{thm:basis}'s hypotheses, hence the theorem applies. 

Then the preliminary remark, combined with the fact that for all $k\in[\![\gamma,p]\!]$
$$
\begin{array}{rll}
\dim \mathbb U_k 
&=\dim   \widetilde{\mathbb P}_k-\dim\mathbb V_{k} &\text{by the direct sum assumption,}\\
&=\dim \ker\mathcal L_*|_{ \widetilde{\mathbb P}_k}+\text{rk}\mathcal L_*|_{ \widetilde{\mathbb P}_k}-\dim\mathbb V_{k} &\text{by the rank nullity theorem,}\\
&=\dim \ker\mathcal L_*|_{ \widetilde{\mathbb P}_k} &\text{by the bijectvity assumption,}
\end{array}
$$
 indeed shows that Algorithm \ref{algo:qTker} computes a generic element of the quasi-Trefftz space $\mathbb Q\mathbb T_p=\ker(\mathcal D_p)$ given any element $U$ of $\mathbb U =\left(  \bigoplus_{n=0}^{\gamma-1} \widetilde{\mathbb P}_n\right) \oplus\left(\bigoplus_{k=\gamma}^{p} {\mathbb U_{k}}\right)$ as input.
Finally, this concludes the proof since 
$$\dim\mathbb Q\mathbb T_p=\sum_{k=\gamma}^p\dim\mathbb U_k+\sum_{\ell=0}^{\gamma-1} \dim\widetilde{\mathbb P}_{\ell}.$$
Moreover, since $\dim \mathbb U_k =\dim   \widetilde{\mathbb P}_k-\dim\mathbb V_{k} =\dim   \widetilde{\mathbb P}_k-\dim   \widetilde{\mathbb P}_{k-\gamma}$, then
$$\dim\mathbb Q\mathbb T_p=\sum_{k=\gamma}^p\dim   \widetilde{\mathbb P}_k-\sum_{k=\gamma}^p\dim   \widetilde{\mathbb P}_{k-\gamma}+\sum_{\ell=0}^{\gamma-1} \dim\widetilde{\mathbb P}_{\ell}$$
which concludes the proof.
\end{proof}


\section{Partial differential operators}
\label{sec:apps}

The hypothesis of Theorem \ref{thm:qTbasis} on the partial differential operator $\mathcal L$ states that the corresponding operator $\mathcal L_*$ consisting of the highest-order terms of $\mathcal L$ with constant coefficients is surjective, and the implementation of the construction algorithm relies on a right inverse for $\mathcal L_*$.
This section determines which partial differential operators satisfy this hypothesis,
 focusing first on constant-coefficient homogeneous operators 
 and then on
 variable-coefficient partial differential operators. 
The outcome is 
 a complete description of Taylor-based polynomial quasi-Trefftz spaces for scalar differential operators with smooth coefficients, including a procedure to compute quasi-Trefftz bases.

\subsection{General framework}
The section focuses on non-trivial constant-coefficient homogeneous operators $\mathcal L_*$, first proving that such operators defined between spaces of homogeneous polynomials are surjective, then providing an algorithm to compute a right inverse.

\begin{prop}\label{prop:surj*}
Given any $\gamma\in\mathbb N$ and any non-trivial constant-coefficient homogeneous differential operator $\mathcal L_*$ of order $\gamma$, namely any operator of the form
$$
\mathcal L_* = P(\partial) 
\text{ for some }
P=
\sum_{|\mathbf i|=\gamma} \alpha_{\mathbf i}{\mathbf X}^{\mathbf i} \in\widetilde{\mathbb P}_\gamma\backslash \{0_\mathbb P\},
$$
then for any $n\in\mathbb N_0$ the corresponding operator
$$
\begin{array}{rccc}
\mathcal L_*:&\widetilde{\mathbb P}_{n+\gamma}&\to &\widetilde{\mathbb P}_n\\
& \Pi&\mapsto &\mathcal L_*\Pi
\end{array}
$$
is surjective.
\end{prop}
As a side note, the result is trivial if $\gamma=0$.
The proof of this proposition leverages a pairing between spaces of homogeneous polynomials.
\begin{lmm}\label{lmm:innerprod}
Consider the pairing on spaces of homogeneous polynomials in $d$ variables defined as follows:
$$
\forall (a,b)\in(\mathbb N_0)^2,
\forall (P,Q)\in\widetilde{\mathbb P}_a\times\widetilde{\mathbb P}_b,
\langle P,Q\rangle_{a\times b}:= P(\partial)\overline Q|_{\mathbf X=0_{\mathbb C^d}}.
$$
where the complex conjugate of a polynomial $Q$ is the polynomial whose coefficients are the complex conjugate of $Q$'s coefficients.
Then
\begin{itemize}
\item $\forall a\in\mathbb N$, $\langle \cdot,\cdot\rangle_{a\times a}$ defines an inner product on $\widetilde{\mathbb P}_a$, and the canonical basis is orthogonal with respect to this inner product;
\item $\forall (a,b,c)\in(\mathbb N_0)^3$ with $c\geq b$,
$\forall (P,Q,R)\in\widetilde{\mathbb P}_a\times\widetilde{\mathbb P}_b\times\widetilde{\mathbb P}_c$,
$\langle P,Q(\partial)R\rangle_{a\times (c-b)}=\langle P\overline Q,R\rangle_{(a+b)\times c}$.
\end{itemize}
\end{lmm}
\begin{proof}
The key polynomial property for this proof is that for any multi-indices $(\mathbf i,\mathbf j)\in(\mathbb N_0^d)^2$
$$
\partial^{\mathbf i} \mathbf X^{\mathbf j} = \frac{\mathbf j!}{(\mathbf j-\mathbf i)!} \mathbf X^{\mathbf j-\mathbf i} \mathbf 1_{\mathbf i\leq\mathbf j}.
$$

Given $a\in\mathbb N$, then $\langle \cdot,\cdot\rangle_{a\times a}$ is a sesquilinear form by definition.  Moreover, for any $(\mathbf i,\mathbf j)\in (\mathbb N_0^d)^2$ with $|\mathbf i|=a$ and $|\mathbf j|=a$, 
\begin{equation}\label{eq:orth}
\left\langle \mathbf X^{\mathbf i}, \mathbf X^{\mathbf j}\right\rangle_{a\times a} 
=\mathbf i! \mathbf 1_{\mathbf i=\mathbf j}
\end{equation}
As a consequence, given $P\in\widetilde{\mathbb P}_a$ defined in the canonical basis by its coefficients $\left\{{p_{\mathbf j}\in\mathbb C,|\mathbf j|=n} \right\}$
$$
\langle P,P\rangle_{a\times a}
=\sum_{|\mathbf i|=n} p_{\mathbf i} \overline{p_{\mathbf i}} 
\geq 0,
$$
where $\langle P,P\rangle_{a\times a}=0$ holds if an only if $P$ is the polynomial 0.
Hence $\langle \cdot,\cdot\rangle_{a\times a}$ is positive definite.
This proves the first point, as orthogonality follows from \eqref{eq:orth}.

Given $(a,b,c)\in(\mathbb N_0)^3$ with $c\geq b$ and $(P,Q,R)\in\widetilde{\mathbb P}_a\times\widetilde{\mathbb P}_b\times\widetilde{\mathbb P}_c$ ,
then $Q(\partial)R\in \widetilde{\mathbb P}_{c-b}$ and 
$$
P(\partial) \overline{\left[ Q(\partial) R\right]}
=
[P\overline{Q}](\partial) \overline{R},
$$
which concludes the proof.
\end{proof}

\begin{proof}[Proof of Theorem \ref{prop:surj*}.]
Consider $n\in\mathbb N$ and $\mathcal L_* = P(\partial) $ for some $P\in\widetilde{\mathbb P}_\gamma$.
Assuming that the operator is not surjective, then there exist $\Lambda \in\widetilde{\mathbb P}_n\backslash \{0_{\mathbb P}\}$ such that 
$$
\forall \Pi\in\widetilde{\mathbb P}_{n+\gamma},
\left\langle  \Lambda,P(\partial)\Pi\right\rangle_{n\times n} =0,
$$
or equivalently, from Lemma \ref{lmm:innerprod}, such that
$$
\forall \Pi\in\widetilde{\mathbb P}_{n+\gamma},
\left\langle \Lambda \overline P, \Pi\right\rangle_{(n+\gamma)\times (n+\gamma)} =0.
$$
In turn, again from Lemma \ref{lmm:innerprod}, this implies that $\Lambda \overline P=0$. This is a contradiction since neither $\Lambda$ nor $P$ can be $0_{\mathbb P}$. This concludes the proof.
\end{proof}
As for the construction of a right inverse operator, 
the task can be performed for a given operator $\mathcal L_*$ depending on the values of its non-zero coefficients.
Consider a non-trivial constant-coefficient homogeneous differential operator  $\mathcal L_*=\sum_{\mathbf i\in I} \alpha_{\mathbf i} \partial^{\mathbf i}$ defined by the sets of its coefficients $\{\alpha_{\mathbf i}\in\mathbb C, |\mathbf i|=\gamma\}$.
Then consider the set of multi-indices of non-zero coefficients $I:=\{\mathbf i\in(\mathbb N_0)^d, |\mathbf i|=\gamma, \alpha_{\mathbf i}\neq 0\}$.
The general idea relies on the definition of
(1) an index $\mathbf i_*\in I$, 
and (2) a linear order $\prec$ on $(\mathbb N_0)^d$
such that
\begin{equation}\label{eq:i*}
 \forall \mathbf i\in I\backslash\{\mathbf i_* \}, \mathbf i\prec\mathbf i_*.
 \end{equation}
For convenience, define the numbering associated to $\prec$, $\mathcal N_\prec:(\mathbb N_0)^d\to\mathbb N$, by:
$$
\forall\mathbf i\in(\mathbb N_0)^d, \mathcal N_\prec(\mathbf i) = \dim\{\mathbf j\in(\mathbb N_0)^d, \mathbf j\prec\mathbf i\}.
$$
The reasoning, leveraging the fact that $\alpha_{\mathbf i_*}\neq 0$, is the following.
For any $n\in\mathbb N_0$ consider the subspace of $\widetilde{\mathbb P}_{n+\gamma}$ defined by
$$
\mathbb V_{n+\gamma}:=\left\{ \Pi\in\widetilde{\mathbb P}_{n+\gamma} ,\Pi=\sum_{|\mathbf j|=n} \pi_{\mathbf j}{\mathbf X}^{\mathbf j+\mathbf i_*}\right\}
\text{ with }
\dim\mathbb V_{n+\gamma} = \dim\widetilde{\mathbb P}_{n};
$$
for any $\Pi\in\mathbb V_{n+\gamma}$ defined by its coefficients in the canonical basis $\left\{{\pi_{\mathbf j}\in\mathbb C,|\mathbf j|=n} \right\}$, then 
$$
\mathcal L_*(\Pi) 
=
\sum_{\mathbf i\in I} 
\sum_{|\mathbf j|=n}
\left(\prod_{l=1}^d\mathbf 1_{i_{l}\leq j_{l}+(i_*)_l}\right)
\alpha_{\mathbf i} \pi_{\mathbf j}\frac{(\mathbf j+\mathbf i_*)!}{(\mathbf j+\mathbf i_*-\mathbf i)!}{\mathbf X}^{\mathbf j+\mathbf i_*-\mathbf i}.
$$
Then, for any $\Lambda\in\widetilde{\mathbb P}_{n}$ defined by its coefficients in the canonical basis $\left\{{\lambda_{\mathbf j}\in\mathbb C,|\mathbf j|=n} \right\}$, $\mathcal L_*\Pi = \Lambda$ is equivalent to
$$
\sum_{|\mathbf j|=n}
\alpha_{\mathbf i_*} \pi_{\mathbf j}\frac{(\mathbf j+\mathbf i_*)!}{\mathbf j!}{\mathbf X}^{\mathbf j}
=
\sum_{|\mathbf j|=n}\lambda_{\mathbf j} \mathbf X^{\mathbf j} 
-
\sum_{\mathbf i\in I\backslash\{\mathbf i_*\}} 
\sum_{|\mathbf j|=n}
\left(\prod_{l=1}^d\mathbf 1_{i_{l}\leq j_{l}+(i_*)_l}\right)
\alpha_{\mathbf i} \pi_{\mathbf j}\frac{(\mathbf j+\mathbf i_*)!}{(\mathbf j+\mathbf i_*-\mathbf i)!}{\mathbf X}^{\mathbf j+\mathbf i_*-\mathbf i},
$$
or equivalently
\begin{equation}\label{eq:itbis}
\forall \mathbf s\in \mathbb N_0^d, |\mathbf s|=n,
\alpha_{\mathbf i_*} \pi_{\mathbf s}\frac{(\mathbf s+\mathbf i_*)!}{\mathbf s!}{\mathbf X}^{\mathbf s}
=
\lambda_{\mathbf s} \mathbf X^{\mathbf s} 
-
\sum_{\mathbf i\in I\backslash\{\mathbf i_*\}} 
\sum_{\mathbf j = \mathbf s-\mathbf i_*+\mathbf i}
\left(\prod_{l=1}^d\mathbf 1_{i_{l}\leq j_{l}+(i_*)_l}\right)
\alpha_{\mathbf i} \pi_{\mathbf j}\frac{(\mathbf j+\mathbf i_*)!}{(\mathbf j+\mathbf i_*-\mathbf i)!}{\mathbf X}^{\mathbf j+\mathbf i_*-\mathbf i}.
\end{equation}
In the right hand side, since $\mathbf i\prec\mathbf i_*$ then $ \mathbf s+\mathbf i-\mathbf i_*\prec \mathbf s$. Therefore, for any given $\mathbf s$, $\pi_{\mathbf s}$ can be computed if $\{\pi_{\mathbf j}, |\mathbf j|=n,\mathbf j\prec \mathbf s\}$ are known.
 As a consequence, given any $\Lambda\in\widetilde{\mathbb P}_{n}$, a solution $\Pi\in\mathbb V_{n+\gamma}$ to $\mathcal L_*\Pi = \Lambda$ can be computed from \eqref{eq:itbis} in the order corresponding to $\mathcal N_\prec$.
For the record, this evidences an underlying triangular structure of the operator $\mathcal L_*$ when described in bases (respectively of $\mathbb V_{n+\gamma}$ and $\widetilde{\mathbb P}_n$) chosen so that elements are ordered according to $\mathcal N_\prec$.
This shows that $\mathcal L_*|_{\mathbb V_{n+\gamma}}$ is surjective,
 and for any $n\in[\![0,p-\gamma]\!]$ it defines a right inverse $\mathcal S_{n+\gamma}$ 
satisfying the hypothesis of Theorem \ref{thm:qTbasis}. The implementation of this right inverse is formalized  in Algorithm \ref{algo:L*gen}.
\begin{algorithm}[!ht]
\label{algo:L*gen}
\DontPrintSemicolon
	{\bf Input} \\
	$\qquad$$\gamma\in\mathbb N$, \\
	$\qquad$ $n\in\mathbb N_0$, \\
	$\qquad$ $k\in [\![1,d]\!]$, \\
	$\qquad$$ \{\alpha_\mathbf i\in\mathbb C,\mathbf i\in\mathbb N_0^d, |\mathbf i| = \gamma, \alpha_{\mathbf i_*}\neq 0\}$ \\
	$\qquad$ $\Lambda\in\widetilde{\mathbb P}_n$ defined by $\left\{{\lambda_{\mathbf j}\in\mathbb C,|\mathbf j|=n} \right\}$
	
        
    \For{$\ell$ from $1$ to $\dim\widetilde{\mathbb P}_{n+\gamma}-\dim \widetilde{\mathbb P}_n$\,}{
    
   $ \mathbf s\gets \mathcal N_\prec^{-1}(\ell)$
    
	    $\displaystyle R_{\mathbf s} \gets\sum_{\mathbf i\in I\backslash\{\mathbf i_*\}} 
\sum_{\mathbf j = \mathbf s-\mathbf i_*+\mathbf i}
\left(\prod_{l=1}^d\mathbf 1_{i_{l}\leq j_{l}+(i_*)_l}\right)
\alpha_{\mathbf i} \pi_{\mathbf j}\frac{(\mathbf j+\mathbf i_*)!}{(\mathbf j+\mathbf i_*-\mathbf i)!}
$

	    	$\displaystyle \pi_{\mathbf s} \gets \frac{\mathbf s!}{(\mathbf s+\mathbf i_*)! \alpha_{\mathbf i_*} }\left(  \lambda_{\mathbf s}
	    	-R_{\mathbf s}\right)$
    }
    
    {\bf Output} $\displaystyle \Pi=\sum_{|\mathbf j|=n} \pi_{\mathbf j}{\mathbf X}^{\mathbf j+\gamma\mathbf e_k}$
    
\caption{Computing a right inverse for $\mathcal L_*$ given $\mathbf i_*$ and $\prec$ satisfying \eqref{eq:i*}}
\end{algorithm}

In practice, one can distinguish different cases. 
If $I=\{\mathbf i\in(\mathbb N_0)^d, |\mathbf i|=\gamma, \alpha_{\mathbf i}\neq 0\}$ contains a single index, or equivalently $\mathcal L_*$ has a single term, the case is trivial: choose $\mathbf i_0$ as the unique index in $I$, then $\mathcal L_*\Pi = \Lambda$ is equivalent to the diagonal system
$$
\sum_{|\mathbf j|=n}
\alpha_{\mathbf i_*} \pi_{\mathbf j}\frac{(\mathbf j+\mathbf i_*)!}{\mathbf j!}{\mathbf X}^{\mathbf j}
=
\sum_{|\mathbf j|=n}\lambda_{\mathbf j} \mathbf X^{\mathbf j} ,
$$
so, given any $\Lambda\in\widetilde{\mathbb P}_{n}$, $\{\pi_{\mathbf j}, |\mathbf j|=n\}$ can be computed in any order. 

One particular case of interest covers a wide number of applications.
If there exists $(k,\mathbf i_*)\in[\![1,d]\!]\times I$, with $\mathbf i_*=(w_1,\cdots,w_d)$, such that for all $\mathbf i\in I\backslash\{\mathbf i_0\}$ then $i_{k}<w_{k}$, then
$$
\mathcal L_*(\Pi) 
= 
\sum_{s=0}^n
\sum_{r=\max(w_{k}-s,0)}^{w_{k} }
\sum_{|\mathbf i|=\gamma} 
\sum_{|\mathbf j|=n}
\mathbf 1_{j_{k}=s-w_{k}+r}\mathbf 1_{i_{k}=r} 
\left(\prod_{l=1}^d\mathbf 1_{i_{l}\leq j_{l}+w_l}\right)
\alpha_{\mathbf i} \pi_{\mathbf j}\frac{(\mathbf j+\mathbf i_*)!}{(\mathbf j+\mathbf i_*-\mathbf i)!}{\mathbf X}^{\mathbf j+\mathbf i_*-\mathbf i},
$$
and, accordingly, $\mathcal L_*\Pi = \Lambda$ is equivalent to
\begin{align*}
&\forall s\in[\![0,n]\!],\\
&
\sum_{|\mathbf j|=n}
\mathbf 1_{j_{k}=s}
\alpha_{\mathbf i_*} \pi_{\mathbf j}\frac{(\mathbf j+\mathbf i_*)!}{\mathbf j!}{\mathbf X}^{\mathbf j}
=
\sum_{|\mathbf j|=n}\lambda_{\mathbf j} \mathbf X^{\mathbf j} \\&
-
\sum_{r=\max(w_{k}-s,0)}^{w_{k}-1 }
\sum_{|\mathbf i|=\gamma} 
\sum_{|\mathbf j|=n}
\mathbf 1_{j_{k}=s-w_{k}+r}\mathbf 1_{i_{k}=r} 
\left(\prod_{l=1}^d\mathbf 1_{i_{l}\leq j_{l}+w_l}\right)
\alpha_{\mathbf i} \pi_{\mathbf j}\frac{(\mathbf j+\mathbf i_*)!}{(\mathbf j+\mathbf i_*-\mathbf i)!}{\mathbf X}^{\mathbf j+\mathbf i_*-\mathbf i}
\end{align*}
Hence given any $\Lambda\in\widetilde{\mathbb P}_{n}$, $\{\pi_{\mathbf s}, |\mathbf s|=n\}$ can be computed for increasing values of the $k$th component of $\mathbf s$. As a side note, this evidences a block triangular structure of $\mathcal L_*$ with diagonal diagonal blocks.
Since this is a particularly common case (more specifically for $\mathbf i_*=\gamma \mathbf e_k$), this is formalized in Algorithm \ref{algo:L*invCC}, a particular version of the general Algorithm \ref{algo:L*gen}.
Note that this includes the trivial case when $I\backslash\{\mathbf i_*\}=\emptyset$.

\begin{algorithm}[!ht]
\label{algo:L*invCC}
\DontPrintSemicolon
	{\bf Input} \\
	$\qquad$$\gamma\in\mathbb N$, \\
	$\qquad$ $n\in\mathbb N_0$, \\
	$\qquad$ $k\in [\![1,d]\!]$, \\
	$\qquad$$ \{\alpha_\mathbf i\in\mathbb C,\mathbf i\in\mathbb N_0^d, |\mathbf i| = \gamma, \alpha_{\mathbf i_*}\neq 0\}$ \\
	$\qquad$ $\Lambda\in\widetilde{\mathbb P}_n$ defined by $\left\{{\lambda_{\mathbf j}\in\mathbb C,|\mathbf j|=n} \right\}$
	
        
    \For{$s$ from $0$ to $n$\,}{
	    \For{$\mathbf s\in\mathbb N_0^d, |\mathbf s| = n, s_k=s$\,}{
	    $\displaystyle R_{\mathbf s} \gets 0$\\
	    \If{$s\geq 1$}{
	    $\displaystyle R_{\mathbf s} \gets \sum_{r=\max(w_k-s,0)}^{w_k-1} \sum_{|\mathbf i|=\gamma} \sum_{|\mathbf j|=n}
	    	\mathbf 1_{\mathbf j=\mathbf s +\mathbf i-\mathbf i_*}
	    	\mathbf 1_{i_k=r} 
	    	\alpha_{\mathbf i} \pi_{\mathbf j}\frac{(\mathbf j+\mathbf i_*)!}{(\mathbf j+\mathbf i_*-\mathbf i)!}$
	    }
	    	$\displaystyle \pi_{\mathbf s} \gets \frac{\mathbf s!}{(\mathbf s+\mathbf i_*)! \alpha_{\mathbf i_*}} \left(  \lambda_{\mathbf s}
	    	-R_{\mathbf s}\right)$
	    }
    }
    
    {\bf Output} $\displaystyle \Pi=\sum_{|\mathbf j|=n} \pi_{\mathbf j}{\mathbf X}^{\mathbf j+\gamma\mathbf e_k}$
    
\caption{Computing a right inverse for $\mathcal L_*$ with $(k,\mathbf i_*)\in[\![1,d]\!]\times I$, with $\mathbf i_*=(w_1,\cdots,w_d)$, such that for all $\mathbf i\in I\backslash\{\mathbf i_0\}$ then $i_{k}<w_{k}$}
\end{algorithm}

Consider now any non-trivial operator $\mathcal L_*$.
In this case, define successively
  $$w_1=\max\{i_1,\mathbf i\in I\},\ I^1:=\{\mathbf i\in I, i_{1}=w_1\},\text{ and }C^1:=\{\mathbf i\in I, i_{1}<w_1\},$$ 
and for $l$ from 2 to $d-1$
 $$w_{l}=\max \{i_l,\mathbf i\in I^{l-1}\},\
I^l:=\{\mathbf i\in I^{l-1}, i_{l}=w_l\},
\text{ and }C^l:=\{\mathbf i\in I^{l-1}, i_{l}<w_l\}.$$
Then $I^l\subset I^{l-1}$ for all  $l$ from $2$ to $d-1$, and  $I^{d-1} = \{\mathbf i\in I,i_{l}=w_l\ \forall l\in [\![1,d-1]\!]\}$ contains a single element, since $|\mathbf i|=\gamma$ implies that  $i_d:=\gamma-\sum_{l=1}^{d-1} i_l$. 
Choose $\mathbf i_*$ to be this single element, with $\mathbf i_*=(w_1,\cdots,w_d)$, so that $I\backslash\{\mathbf i_*\} = \cup_{l=1}^{d-1} C^l$.
Note that the case $\cup_{l=2}^{d-1} C^l=\emptyset$ is the common case described previously.
Considering
the linear order $\prec$ on $ \mathbb N_0^d$,  defined by	
$$
\forall (\mu,\nu)\in( \mathbb N_0^d)^2,
 \mu\prec\nu 
 \Leftrightarrow
\left\{
 \begin{array}{l}
 \mu_1<\nu_1\ , \text{ or }\\
\exists k\in[\![2,d]\!],  \forall k'\in[\![1,k-1]\!],  \mu_{k'}=\nu_{k'} \text{ and } \mu_{k}<\nu_{k}.
 \end{array}
 \right.
$$
then $\forall \mathbf i\in I\backslash\{\mathbf i_* \}, \mathbf i\prec\mathbf i_*$, the desired property \eqref{eq:i*}.
Algorithm \ref{algo:L*gen} can then be applied for these choices of $\mathbf i_*$ and $\prec$.

\subsection{Quasi-Trefftz framework}
Combining the different results from this work leads to the following theorem which holds for any non-trivial scalar linear differential operator with smooth coefficients.

\begin{thm}\label{thm:qTccl}
Given $(d,\gamma,p)\in\mathbb N^3$ such that $p>\gamma$ and $\mathbf x_0\in\mathbb R^d$,
consider
a linear differential operator with variable coefficients $\mathbf c = \{c_\mathbf j\in\mathcal C^{p-\gamma}(\mathbf x_0),\mathbf j\in\mathbb N_0^d, |\mathbf j| \leq \gamma\}$, $\displaystyle \mathcal L := \sum_{m=0}^\gamma \sum_{ |\mathbf j|=m} c_{\mathbf j} 
 \partial_{{\mathbf x}}^{\mathbf j}$ and define $\mathcal L_*: = \sum_{ |\mathbf j|=\gamma} c_{\mathbf j}\left( {\mathbf x}_0  \right) \partial_{{\mathbf x}}^{\mathbf j}$.
 Assume that there exists $\mathbf j\in\mathbb N_0^d$ with $|\mathbf j|=\gamma$ such that $c_\mathbf j (\mathbf x_0)\neq 0$, so that at $\mathbf x_0$ the operator is of order $\gamma$.
 Denote by $I:=\{\mathbf i\in(\mathbb N_0)^d, |\mathbf i|=\gamma, c_\mathbf i (\mathbf x_0)\neq 0\}$ the set of multi-indices of non-zero coefficients of $\mathcal L$ evaluated at $\mathbf x_0$.
 
Consider $\mathbf i_*\in I$ 
and a linear order $\prec$ on $(\mathbb N_0)^d$
such that for all $\mathbf i\in I\backslash\{\mathbf i_* \}$, $\mathbf i\prec\mathbf i_*$. For all $n\in[\![\gamma,p]\!]$ define the polynomial spaces 
$\mathbb V_{n}:=Span\{\mathbf X^{\mathbf i_*+\mathbf i},\mathbf i\in\mathbb N_0^d,|\mathbf i|=n-\gamma\}$ and
$\mathbb U_n$ its complement in $\widetilde{\mathbb P}_n$, i.e. $\widetilde{\mathbb P}_n=\mathbb U_n\oplus \mathbb V_n$, 
and denote by $\mathcal S_{n}$ the right inverse of $\mathcal L_*|_{\mathbb V_n}$ provided by Algorithm \ref{algo:L*gen}.
  
 Then, the quasi-Trefftz operator $\mathcal D_p = T_{p-\gamma}\circ\mathcal L:\mathbb P_p\to\mathbb P_{p-\gamma}$, is surjective, so the dimension of the quasi-Trefftz space $\dim\mathbb Q\mathbb  T_p =\ker \mathcal D_p $ is given by
 $$\dim\mathbb Q\mathbb  T_p = \dim\mathbb P_{p}-\dim\mathbb P_{p-\gamma}.$$
 Moreover, a basis of the quasi-Trefftz space $\mathbb Q\mathbb T_p$ can be computed starting from any basis
 $\{Y_k,k\in [\![1, \dim\mathbb P_{p}-\dim\mathbb P_{p-\gamma}]\!]\}$ of $\mathbb U =\left(  \bigoplus_{n=0}^{\gamma-1} \widetilde{\mathbb P}_n\right) \oplus\left(\bigoplus_{k=\gamma}^{p} \mathbb U_k\right)$,
consider the corresponding  $\{X_k,k\in  [\![1, \dim\mathbb P_{p}-\dim\mathbb P_{p-\gamma}]\!]\}$ where, for each $k$, $X_k$ is the output of Algorithm \ref{algo:qTker} for the input $Y_k$: $\{X_k+Y_k,k\in  [\![1, \dim\mathbb P_{p}-\dim\mathbb P_{p-\gamma}]\!]\}$ forms a basis of $\mathbb Q\mathbb  T_p$.  

 \end{thm}
 \begin{proof}
 Under the hypothesis, then the corresponding constant coefficient operator of order $\gamma$ defined by $\mathcal L_* = \sum_{ |\mathbf j|=\gamma} c_{\mathbf j}\left( {\mathbf x}_0  \right) \partial_{{\mathbf x}}^{\mathbf j}$ is non trivial. As a consequence, for all $n\in\mathbb N_0$, $\mathcal L_* |_{\widetilde{\mathbb P}_{n+\gamma}}:\widetilde{\mathbb P}_{n+\gamma}\to\widetilde{\mathbb P}_n$ is surjective according to Proposition \ref{prop:surj*}, and more precisely $\mathcal L_* |_{{\mathbb V}_{n+\gamma}}$ is bijective and $\mathcal S_n$ is a right inverse.
 
 Since for all $n\in[\![\gamma,p]\!]$ $\mathbb V_{n}=Span\{\mathbf X^{\mathbf i},\mathbf i\in\mathbb N_0^d,|\mathbf i|=n, \mathbf i\geq\mathbf i_* \}$,
it is clear that $\widetilde{\mathbb P}_n=\mathbb U_n\oplus \mathbb V_n$. The conclusions then simply follow from combining Theorem \ref{thm:qTbasis} and Proposition \ref{prop:surj*}.
 \end{proof}
 
Generally, second order equations of order exactly 2 at a point $\mathbf x_0$ can be written as
$$\nabla\cdot (\mathbf A(\mathbf x) \nabla u(\mathbf x)+ \mathbf b(\mathbf x) u(\mathbf x))+n(\mathbf x)u(\mathbf x)=0,$$
 where $\mathbf A$ is matrix-valued and at least one of its entries is non-zero at $\mathbf x_0$, and $\mathbf b$ is vector-valued.
A wide range of standard second order PDEs, including elliptic equations, parabolic equations, hyperbolic equations, and type-changing equations fall under particular case for $\mathcal L_*$, i.e. there exists $(k,\mathbf i_*)\in[\![1,d]\!]\times I$, with $\mathbf i_*=(w_1,\cdots,w_d)$, such that for all $\mathbf i\in I\backslash\{\mathbf i_0\}$ then $i_{k}<w_{k}$. For these equations, Algorithm \ref{algo:L*gen} can be replaced by the simplest Algorithm \ref{algo:L*invCC}.
As an illustration, this is the case of:
\begin{itemize}
\item the anisotropic Laplace equation $-\nabla\cdot( \mathbf A(\mathbf x) \nabla u(\mathbf x))=0$, where $\mathbf A$ is matrix-valued and has all eigenvalues of the same sign at $\mathbf x_0$,
\item the anisotropic Helmholtz equation $-\nabla\cdot( \mathbf A(\mathbf x) \nabla u(\mathbf x))-\kappa^2n(\mathbf x)u(\mathbf x)=0$, where $\mathbf A$ is matrix-valued and has positive eigenvalues at $\mathbf x_0$,
\item the convected Helmholtz equation $-\nabla\cdot(\rho(\mathbf x) (\nabla u(\mathbf x) -(\mathbf M(\mathbf x) \cdot\nabla u(\mathbf x) )\mathbf M(\mathbf x) +\mi \kappa u(\mathbf x) \mathbf M(\mathbf x) )
-\rho(\mathbf x) (\kappa^2u(\mathbf x) +\mi \kappa \mathbf M(\mathbf x) \cdot\nabla u(\mathbf x) )=0$, where $\kappa>0$,  and at $\mathbf x_0$: $\mathbf M$ is vector-valued, $\|\mathbf M\|<1$, and $\rho$ is positive,
\item the diffusion-advection-reaction equation $-\nabla\cdot (\mathbf A(\mathbf x) \nabla u(\mathbf x)+ \mathbf b(\mathbf x) u(\mathbf x))+n(\mathbf x)u(\mathbf x)=0$, where $\mathbf A$ is matrix-valued and has positive eigenvalues at $\mathbf x_0$, and $\mathbf b$ is vector-valued,
\item the wave equation $\partial_t^2u(\mathbf x,t)-\nabla\cdot( \mathbf A(\mathbf x,t) \nabla u(\mathbf x,t))=0$, where $\mathbf A$ is matrix-valued and  and has positive eigenvalues at $\mathbf x_0$,
\item the heat equation $\partial_tu(\mathbf x,t)-\nabla\cdot( \mathbf A(\mathbf x,t) \nabla u(\mathbf x,t))=0$, where $\mathbf A$ is matrix-valued and  and has positive eigenvalues at $\mathbf x_0$,
\item the advection equation $\partial_tu(\mathbf x,t)+ \mathbf b(\mathbf x)\cdot \nabla u(\mathbf x,t)=0$, where $\mathbf b$ is vector-valued,
\item the Tricomi equation $\partial_x^2 u(x,y)+x\partial_y^2 u(x,y)=0$ at any point either in the hyperbolic half plane$ x < 0$, in the elliptic half plane $x > 0$, or on the degenerates line $x = 0$,
\item the Schr\"odinger equation with variable potential $\mathrm i\partial_tu(\mathbf x,t)+\frac 12\Delta_{\mathbf x} u(\mathbf x,t)-V(\mathbf x,t) u(\mathbf x,t)=0$,
\end{itemize}

\begin{rmk}
Consider a linear differential operator of order at most equal to $\gamma$, $\displaystyle \mathcal L := \sum_{m=0}^{\gamma} \sum_{ |\mathbf j|=m} c_{\mathbf j}  \partial_{{\mathbf x}}^{\mathbf j}$, defined by its variable coefficients $\mathbf c = \{c_\mathbf j\in\mathcal C^{p-\gamma},\mathbf j\in\mathbb N_0^d, |\mathbf j| \leq \gamma\}$.
If at some point $\mathbf x_0$ there exists $(\tilde\gamma,\mathbf j)\in\mathbb N\times\mathbb N_0^d$ with $\tilde\gamma<\gamma$, $|\mathbf j|=\tilde\gamma$ and $c_\mathbf j (\mathbf x_0)\neq 0$ while $c_\mathbf j (\mathbf x_0)= 0$ for all $\mathbf j$ such that $\tilde\gamma<|\mathbf j|\leq \gamma$, then 
the Theorem can be applied to $\displaystyle \tilde{ \mathcal L} := \sum_{m=0}^{\tilde\gamma} \sum_{ |\mathbf j|=m} c_{\mathbf j}  \partial_{{\mathbf x}}^{\mathbf j}$. So in this situation $\dim\mathbb Q\mathbb  T_p = \dim\mathbb P_{p}-\dim\mathbb P_{p-\tilde\gamma}$ and the theorem also provides an explicit procedure to construct a quasi-Trefftz basis.
\end{rmk}

\section{Comments and conclusion}
\label{sec:ccl}
This article focuses on local Taylor-based polynomial quasi-Trefftz spaces for homogeneous scalar linear equations with smooth coefficients, and proposes an explicit procedure to compute quasi-Trefftz basis.
Note that a minimal adjustment to this proposed procedure yields a procedure to compute a particular quasi-Trefftz solution to a scalar linear equations with smooth coefficients and a smooth right hand side.

For a homogeneous scalar linear equations with smooth coefficients defined over a domain in $\mathbb R^d$, global quasi-Trefftz spaces can be constructed over the domain by simply meshing the domain and choosing one point $\mathbf x_0$ per mesh element to construct quasi-Trefftz bases with support restricted to each element. 
Obtaining global approximation properties from local approximation properties requires assumptions on the mesh, see for instance \cite{IGMPS-dar}.


Beyond the framework of polynomial quasi-Trefftz spaces for scalar equations, the general framework developed in this article lays the ground to tackle two significant challenges.
The first one concerns the use of different ansatz to define quasi-Trefftz functions. An ansatz tailored for a particular governing PDE might . This is illustrated for the case of GPWs in \cite{GPWgen}.
The second one concerns the development of some quasi-Trefftz methods for systems of PDEs. For comparison, Trefftz methods have been studied for Maxwell problems \cite{Seb,MoiolaMax} or for Friedrichs systems \cite{morelBFs}, with bases of plane waves, of polynomials or of exponentials. By contrast, the main difficulty in systematically constructing quasi-Trefftz bases for systems of PDEs lies in the fact that homogeneous vector-valued differential operators are in general surjective, think for instance of the gradient operator. In this case, vector-valued spaces of polynomials have a graded structure to scalar-valued ones. Hence the general framework developed here is a natural starting point.

\section{Aknowledgement}
I would like to acknowledge support from the US National Science Foundation (NSF):
this material is based upon work supported by the NSF under Grant No. DMS-2110407.

I would also like to warmly thank Bruno Despr\'es, Anton Izosimov and Tonatiuh S\'anchez Vizuet for fruitful conversations that helped me shape the material presented in this article.

\section*{Statements and declarations}
 L.-M. Imbert-G\'erard has disclosed an outside interest in Airbus Central R\&T to the University of Arizona.  Conflicts of interest resulting from this interest are being managed by The University of Arizona in accordance with its policies.

%

\appendix
\section{Some useful formulas}

\subsection{Explicit formula for the quasi-Trefftz operator}
\label{app:newexpl}
The definition of the quasi-Trefftz operator $\mathcal D  = T_{p-\gamma}\circ \mathcal L$ can be expressed explicitly as follows.

Given the variable coefficients of $\mathcal L$ $\{ c_{
n}\in\widetilde{\mathbb P}_{n},n\in [\![0,p-\gamma]\!] \}$, define $c_{\mathbf j}=\sum_{n=0}^{p-\gamma}c_{n}$, then for any $\Pi_N\in\widetilde{\mathbb P}_N$ 
$$
c_{\mathbf j}\partial_{\mathbf x}^{\mathbf j}\Pi_N = 
\left(\sum_{n=0}^{p-\gamma}c_{n}\right)\partial_{\mathbf x}^{\mathbf j}\Pi_N = 
\left\{\begin{array}{ll}
0_{\mathbb P} & \text{ if } |\mathbf j|>N,\\
\displaystyle\sum_{\ell = N-|\mathbf j|}^{p-\gamma+N-|\mathbf j|} c_{\ell-(N-|\mathbf j|)} \partial_{\mathbf x}^{\mathbf j}\Pi_N & \text{ if } |\mathbf j|\leq N,
\end{array}\right.
$$
where $c_{\ell-(N-|\mathbf j|)}\in\widetilde{\mathbb P}_{\ell-(N-|\mathbf j|)}$ and $\partial_{\mathbf x}^{\mathbf j}\Pi_N \in\widetilde{\mathbb P}_{N-|\mathbf j|} $,  so $ c_{\ell-(N-|\mathbf j|)} \partial_{\mathbf x}^{\mathbf j}\Pi_N \in\widetilde{\mathbb P}_{\ell}$.
Hence
$$
T_{p-\gamma} \left(c_{\mathbf j}\partial_{\mathbf x}^{\mathbf j}\Pi_N \right)= 
\left\{\begin{array}{ll}
0_{\mathbb P} & \text{ if } |\mathbf j|>N,\\
0_{\mathbb P} & \text{ if } N-|\mathbf j|>p-\gamma,\\
\displaystyle\sum_{\ell = N-|\mathbf j|}^{p-\gamma} c_{\ell-(N-|\mathbf j|)} \partial_{\mathbf x}^{\mathbf j}\Pi_N & \text{ if }N-(p-\gamma)\leq  |\mathbf j|\leq N.
\end{array}\right.
$$
As a result, given  $\{ c_{\mathbf j,n}\in\widetilde{\mathbb P}_{n},n\in [\![0,p-\gamma]\!] , \mathbf j \in (\mathbb N_0)^d, |\mathbf j|\leq \gamma\}$, define $c_{\mathbf j}=\sum_{n=0}^{p-\gamma}c_{\mathbf j,n}$, 
if $\displaystyle \mathcal L := \sum_{m=0}^\gamma \sum_{ |\mathbf j|=m} c_{\mathbf j}\left( {\mathbf x}  \right) \partial_{{\mathbf x}}^{\mathbf j}$ then for any $\Pi_N\in\widetilde{\mathbb P}_N$ 
\begin{equation}\label{eq:dpiN}
\mathcal D(\Pi_N) = \sum_{m=\max(0,N-(p-\gamma))}^{\min(N,\gamma)} \sum_{ |\mathbf j|=m}  \sum_{\ell=N-m}^{p-\gamma} c_{\mathbf j,\ell-(N-m)} \partial_{\mathbf x}^{\mathbf j}\Pi_N
\end{equation}
where here again $ c_{\mathbf j,\ell-(N-|\mathbf j|)} \partial_{\mathbf x}^{\mathbf j}\Pi_N \in\widetilde{\mathbb P}_{\ell}$. 
In particular, 
\begin{itemize}
\item if $N\geq \gamma$, then the smallest value of $\ell$ in the last sum is $\ell = N-\gamma$, so $\mathcal D(\Pi_N)\in \bigoplus_{\ell = N-\gamma}^{p-\gamma}\widetilde{\mathbb P}_\ell$;
\item if $N<\gamma$, then the smallest value of $\ell$ in the last sum is $\ell = 0$, so $\mathcal D(\Pi_N)\in \bigoplus_{\ell = 0}^{p-\gamma}\widetilde{\mathbb P}_\ell$
\end{itemize}

\subsection{Additive decomposition of the operator}
\label{app:newAdd}
Given the operator $\mathcal L_* = \sum_{ |\mathbf j|=\gamma} c_{\mathbf j}\left( {\mathbf x}_0  \right) \partial_{{\mathbf x}}^{\mathbf j}$ and the corresponding remainder operator $\mathcal R = \mathcal D_p-\mathcal L_*$, 
then on the one hand for any $N\in[\![0,p]\!]$ and any $\Pi_N\in\widetilde{\mathbb P}_N$
$$
\mathcal L_*(\Pi_N)=
\left\{\begin{array}{ll}
0_{\mathbb P} & \text{ if } N<\gamma,\\
\displaystyle
 \sum_{ |\mathbf j|=\gamma} c_{\mathbf j,0} \partial_{\mathbf x}^{\mathbf j}\Pi_N& \text{ if } N\geq \gamma, \text{ from } (m,\ell) =( \gamma,N-\gamma) \text{  in \eqref{eq:dpiN},}
\end{array}\right.
$$
while on the other hand  for any $N\in[\![0,p-1]\!]$ and any $\Pi_N\in\widetilde{\mathbb P}_N$
$$
\mathcal R  (\Pi_N) = 
\left\{\begin{array}{rl}
\displaystyle
 \sum_{m=\max(0,N-(p-\gamma))}^{N} \sum_{ |\mathbf j|=m}  \sum_{\ell=N-m}^{p-\gamma} c_{\mathbf j,\ell-(N-m)} \partial_{\mathbf x}^{\mathbf j}\Pi_N & \text{ if } N<\gamma,
\\
\displaystyle
\sum_{m=\max(0,N-(p-\gamma))}^{\gamma-1} \sum_{ |\mathbf j|=m}  \sum_{\ell=N-m}^{p-\gamma} c_{\mathbf j,\ell-(N-m)} \partial_{\mathbf x}^{\mathbf j}\Pi_N\\
\displaystyle
+\sum_{ |\mathbf j|=\gamma}  \sum_{\ell=N-\gamma+1}^{p-\gamma} c_{\mathbf j,\ell-(N-m)} \partial_{\mathbf x}^{\mathbf j}\Pi_N & \text{ if } N\geq \gamma,
\end{array}\right.
$$
and finally for $N=p$ then \eqref{eq:dpiN} reads $\mathcal D_p(\Pi_p) = \mathcal L_*(\Pi_p)$ which implies that $\mathcal R(\Pi_p) = 0_{\mathbb P_{p-\gamma}}$.

\end{document}